\documentclass[a4paper,11pt]{article}

\usepackage{amsmath}
\usepackage{amssymb}
\usepackage{amsthm}
\usepackage{color}
\usepackage{cite}
\usepackage{graphicx}
\usepackage{fullpage}
\usepackage{numprint}
\usepackage{url}
\usepackage{subfig}
\usepackage{hyperref}

\npthousandsep{,}\npthousandthpartsep{}\npdecimalsign{.}

\newcommand{\seqk}[2]{\ensuremath{\left( #2 \right)_{#1 \in \N}}}
\newcommand{\set}[2]{\ensuremath{\left\{#1 \; | \; #2\right\}}}
\newcommand{\directset}[1]{\ensuremath{\left\{ #1 \right\}}}

\newcommand{\intrange}[3]{\ensuremath{#1 = #2, \ldots, #3}}
\newcommand{\abs}[1]{\ensuremath{\left| #1 \right|}}
\newcommand{\norm}[1]{\ensuremath{\left\| #1 \right\|}}
\newcommand{\cnorm}[1]{\ensuremath{\norm{#1}_\infty}}

\newcommand{\closure}[1]{\ensuremath{\overline{#1}}}
\newcommand{\convex}[1]{\ensuremath{\mathrm{co}\left(#1\right)}}

\newcommand{\num}[1]{\numprint{#1}}

\newcommand{\B}[2]{\ensuremath{B_{#1}\left(#2\right)}}
\newcommand{\Bc}[2]{\ensuremath{\closure{\B{#1}{#2}}}}
\newcommand{\bigO}[1]{\ensuremath{O\left(#1\right)}}
\newcommand{\floor}[1]{\ensuremath{\left\lfloor #1 \right\rfloor}}
\newcommand{\ceil}[1]{\ensuremath{\left\lceil #1 \right\rceil}}

\newcommand{\integral}[4]{\ensuremath{\int_{#1}^{#2}  #3 \, #4}}
\newcommand{\exptval}[1]{\ensuremath{\mathbb{E}\left(#1\right)}}

\newcommand{\C}[0]{\ensuremath{C}}
\newcommand{\R}[0]{\ensuremath{\mathbb{R}}}
\newcommand{\Q}[0]{\ensuremath{\mathbb{Q}}}
\newcommand{\N}[0]{\ensuremath{\mathbb{N}}}
\renewcommand{\d}[1]{\ensuremath{d_{#1}}}
\newcommand{\sd}[1]{\ensuremath{\mathrm{sd}_{#1}}}
\newcommand{\md}[1]{\ensuremath{\mathrm{md}_{#1}}}
\newcommand{\dH}[0]{\ensuremath{d_H}}
\newcommand{\dHapprox}[0]{\ensuremath{\tilde{d}}}
\newcommand{\dHone}[2]{\ensuremath{d\left(#1 \to #2\right)}}
\newcommand{\gammafcn}[1]{\ensuremath{\Gamma\left(#1\right)}}
\newcommand{\Beta}[2]{\ensuremath{\mathrm{B}\left(#1, \, #2\right)}}

\newenvironment{alternatives}[0]
  {\left\{\begin{array}{rl}}
  {\end{array}\right.}
\newcommand{\alternative}[2]{#1 & #2 \\}

\theoremstyle{definition}
\newtheorem{realDefinition}{Definition}
\newenvironment{definition}[1]
  {\begin{realDefinition}\label{dfn:#1}}
  {\end{realDefinition}}

\newtheorem{realLemma}{Lemma}
\newenvironment{lemma}[1]
  {\begin{realLemma}\label{lemma:#1}}
  {\end{realLemma}}
\newtheorem{realCorollary}{Corollary}
\newenvironment{corollary}[1]
  {\begin{realCorollary}\label{cor:#1}}
  {\end{realCorollary}}
\newtheorem{realThm}{Theorem}
\newenvironment{theorem}[1]
  {\begin{realThm}\label{thm:#1}}
  {\end{realThm}}

\theoremstyle{definition}
\newtheorem{realExample}{Example}
\newenvironment{example}[1]
  {\begin{realExample}
   \label{ex:#1}}
  {\end{realExample}}

\newcommand{\includefig}[4]
  {
    \begin{figure}
    \begin{center}
    \includegraphics[#2]{#1}
    \end{center}
    \caption{#4}
    \label{fig:#3}
    \end{figure}
  }
\newcommand{\figid}[0]{}
\newcommand{\includesubfig}[3]
  {
    \begin{figure}
    \renewcommand{\figid}[0]{#1}
    \begin{center}
      #3
    \end{center}
    \caption{#2}
    \label{fig:#1}
    \end{figure}
  }
\newcommand{\subfig}[4]
  {
    \subfloat[#4]
      {\label{subfig:\figid:#3}\includegraphics[width=#1\textwidth]{#2}}
  }

\newcommand{\mysection}[2]{\section{#2}\label{sec:#1}}
\newcommand{\mysubsection}[2]{\subsection{#2}\label{subsec:#1}}

\renewcommand{\eqref}[1]{\hyperref[eq:#1]{(\ref*{eq:#1})}}
\newcommand{\lemmaref}[1]{\hyperref[lemma:#1]{Lemma~\ref*{lemma:#1}}}
\newcommand{\corref}[1]{\hyperref[cor:#1]{Corollary~\ref*{cor:#1}}}
\newcommand{\thref}[1]{\hyperref[thm:#1]{Theorem~\ref*{thm:#1}}}
\newcommand{\defref}[1]{\hyperref[dfn:#1]{Definition~\ref*{dfn:#1}}}
\newcommand{\exref}[1]{\hyperref[ex:#1]{Example~\ref*{ex:#1}}}
\newcommand{\figref}[1]{\hyperref[fig:#1]{Figure~\ref*{fig:#1}}}
\newcommand{\subfigref}[2]{\hyperref[subfig:#1:#2]{Figure~\ref*{subfig:#1:#2}}}
\newcommand{\secref}[1]{\hyperref[sec:#1]{Section~\ref*{sec:#1}}}
\newcommand{\subsecref}[1]{\hyperref[subsec:#1]{Subsection~\ref*{subsec:#1}}}

\begin{document}

\title{Computing the Hausdorff Distance of Two Sets
from Their Signed Distance Functions}

\author{Daniel Kraft \\
University of Graz \\
Institute of Mathematics, NAWI Graz \\
Universit\"atsplatz 3, 8010 Graz, Austria \\
Email: daniel.kraft@uni-graz.at}

\date{August 4th, 2016}

\maketitle

\begin{abstract}
The Hausdorff distance is a measure of (dis-)similarity between two sets
which is widely used in various applications.  Most of the applied
literature is devoted to the computation for sets consisting of a finite
number of points.  This has applications, for instance, in image processing.
However, we would like to apply the Hausdorff distance to control and evaluate
optimisation methods in level-set based shape optimisation.
In this context, the involved sets are not finite point sets but
characterised by level-set or signed distance functions.
This paper discusses the computation of the Hausdorff distance
between two such sets.  We recall fundamental properties
of the Hausdorff distance, including a characterisation in terms of
distance functions.  In numerical applications, this result gives
at least an exact lower bound on the Hausdorff distance.  We also derive
an upper bound, and consequently a precise error estimate.  By giving
an example, we show that our error estimate cannot be further improved
for a general situation.
On the other hand, we also show that
much better accuracy can be expected for non-pathological
situations that are more likely to occur in practice.
The resulting error estimate can be improved even further
if one assumes that the grid is rotated randomly with respect to the
involved sets.

\vspace{1em}
\end{abstract}

\mysection{introduction}{Introduction}

The Hausdorff distance (also called Pompeiu-Hausdorff distance)
is a classical measure for the difference between two sets:

\begin{definition}{dH}
Let $A, B \subset \R^n$.  Then the \emph{one-sided Hausdorff distance}
between $A$ and $B$ is defined as
\begin{equation}
\label{eq:dHone}
\dHone{A}B = \sup_{x \in A} \inf_{y \in B} \abs{x - y}.
\end{equation}
This allows us to introduce the \emph{Hausdorff distance}:
\begin{equation}
\label{eq:dH}
\dH(A, B) = \max\left(\dHone{A}B, \, \dHone{B}A\right)
\end{equation}
\end{definition}

While one can, in fact, define the Hausdorff distance between
subsets of a general metric space, we are only interested in subsets
of $\R^n$ in the following.
Note that $\dHone{A}B \ne \dHone{B}A$ in general, such that
the additional symmetrisation step in \eqref{dH} is necessary.
For instance, if $A \subset B$, then $\dHone{A}B = 0$ while
$\dH(A, B) = \dHone{B}A > 0$ unless $\closure{A} = \closure{B}$.
Since the Euclidean norm $\abs{\cdot}$ is continuous, it is easy to see that
\eqref{dHone} and thus also $\dH$ is not changed if we replace
one or both of the sets by their interior or closure.
The set of compact subsets of $\R^n$ is turned into a metric space by $\dH$.
For some general discussion about the Hausdorff distance,
see Subsection~6.2.2 of \cite{shapesAndGeometry}.
The main theoretical properties that we need will be discussed
in \secref{properties} in more detail.

Historically, the Hausdorff distance is a relatively old concept.
It was already introduced by Hausdorff in 1914, with a similar concept
for the reciprocal distance between two sets dating back to Pompeiu in 1905.
In recent decades, the Hausdorff distance has found plenty of applications
in various fields.  For instance, it has been applied in image processing
\cite{comparingImagesHausdorff}, object matching \cite{objectMatchingHausdorff},
face detection \cite{faceDetectionHausdorff}
and for evolutionary optimisation \cite{evolutionaryOptHausdorff},
to name just a few selected areas.
In most of these applications, the sets whose distance is computed
are finite point sets.  Those sets may come, for instance,
from filtering a digital image or a related process.
Consequently, there exists a lot of literature that deals with the computation
of the Hausdorff distance for point sets, such as
\cite{incrementalHausdorffAlgorithm}.
Methods exist also for sets of other structure, for instance,
convex polygons \cite{hausdorffConvexPolygons}.

We are specifically interested in applying the Hausdorff distance
to measure and control the progress of level-set based shape optimisation
algorithms such as the methods employed in \cite{breakwaterPaper}
and \cite{mmar2014}.
In particular, the Hausdorff distance between
successive iterates produced by some descent method may be useful to implement
a stopping criterion or to detect when a descent run is getting stuck
in a local minimum.
For these applications, the sets $A$ and $B$ are typically open domains
that are described by the sub-zero level sets of some functions.
To the best of our knowledge, no analysis has been done so far
on the computation of the Hausdorff distance for sets given in this way.
A special choice for the level-set function of a domain is
its \emph{signed distance function}:

\begin{definition}{distance functions}
Let $\Omega \subset \R^n$ be bounded.  We define the
\emph{distance function} of $\Omega$ as
\begin{equation*}
\d\Omega(x) = \inf_{y \in \Omega} \abs{x - y}.
\end{equation*}
Note that $\d\Omega(x) = 0$ for all $x \in \closure\Omega$.
To capture also information about the interior of $\Omega$, we
introduce the \emph{signed} (or oriented) distance function as well:
\begin{equation*}
\sd\Omega(x)
  = \begin{alternatives}
      \alternative{\d\Omega(x)}{x \not\in \Omega,}
      \alternative{-\d{\R^n \setminus \Omega}(x)}{x \in \Omega}
    \end{alternatives}
\end{equation*}
See also chapters 6 and 7 of \cite{shapesAndGeometry}.
Both $\d\Omega$ and $\sd\Omega$ are Lipschitz continuous with constant one.
\end{definition}

If the signed distance function of $\Omega$ is not known,
it can be calculated very efficiently from an arbitrary level-set function
using the Fast Marching Method \cite{sethianFastMarching}.
Conveniently, the Hausdorff distance can be characterised in terms
of distance functions.
It should not come as a big surprise that this is possible,
considering that the distance function $\d{B}(x)$ appears on the right-hand
side of \eqref{dHone}.
This is a classical result, which we will recall
and discuss in \secref{properties}.
For numerical calculations, though, the distance functions are known
only on a finite number of grid points.  In this case, the classical
characterisation only yields an exact \emph{lower bound} for $\dH(A, B)$.
The main result of this paper is the derivation of upper bounds as well,
such that the approximation error can be estimated.
These results will be presented in \secref{error estimates}.
We also give an example to show that our estimates of \subsecref{worst case}
are sharp in the general case.
In addition, we will see that much better estimates can be achieved
for (a little) more specific situations.
Since these situations still cover a wide range of sets that may occur
in practical applications, this result is also useful.
\subsecref{numerics} gives a comparison of the actual numerical
error for some situations in which $\dH(A, B)$ is known exactly.
We will see that these results
match the theoretical conclusions quite well.
In \secref{randomisation}, finally, we show that further improvements
are possible if we assume that the orientation of the grid
is not related to the sets $A$ and $B$.  This can be achieved,
for instance, by a random rotation of the grid, and is usually justified
if the data comes from a measurement process.

Note that our code for the computation of (signed) distance functions
as well as the Hausdorff distance following the method suggested here
has been released as free software.  It is included in the
\texttt{level-set} package \cite{package} for GNU~Octave \cite{octave}.

\mysection{properties}{Characterising $\dH$ in Terms of Distance Functions}

Let $A, B \subset \R^n$ be two compact sets throughout the remainder
of the paper.  In this case, it is easy to see that compactness
implies that the various suprema and infima in \defref{dH}
and \defref{distance functions} are actually attained:


\begin{lemma}{inf/sup attained}
For each $x \in \R^n$ there exist $y_1, y_2 \in A$ such that
$\d{A}(x) = \abs{x - y_1}$ and $\abs{\sd{A}(x)} = \abs{x - y_2}$.
Furthermore, there also exist $x \in A$ and $y \in B$
such that $\dHone{A}B = \abs{x - y}$.
This, of course, implies that $\dH(A, B)$ can also
be expressed in a similar form.
\end{lemma}

Let us now, for the rest of this section, turn our
attention to the relation between the Hausdorff distance
and distance functions.  While most of this content
is well-known and not new, we believe that it makes sense
to give a comprehensive discussion.  This is particularly
true because the Hausdorff distance is a concept that can be
quite unintuitive.  Thus, we try to clearly explain
potential pitfalls and give counterexamples where appropriate.
This discussion forms the basis for the later sections, in which we
present our new results.

\mysubsection{dH and dist}{Distance Functions}

One may have the idea to ``characterise'' the sets $A$ and $B$
via their distance functions $\d{A}$ and $\d{B}$ from
\defref{distance functions}.
Since the distance functions are part of the Banach space
$\C(\R^n)$ of continuous functions, the norm on this space
can be used to define a distance between $A$ and $B$ as
$\cnorm{\d{A} - \d{B}}$.
We will now see that this distance is equal to the
Hausdorff distance defined in \defref{dH}:

\begin{theorem}{dH as sup over d}
For each $x \in \R^n$, the inequality
\begin{equation}
\label{eq:dH as sup over d inequality}
\abs{\d{A}(x) - \d{B}(x)} \le \dH(A, B)
\end{equation}
holds.  More precisely, one even has
\begin{equation}
\label{eq:dH as sup over d}
\dH(A, B) = \cnorm{\d{A} - \d{B}}
  = \sup_{x \in \R^n} \abs{\d{A}(x) - \d{B}(x)}.
\end{equation}

\begin{proof}
This is a classical result, which is, for instance,
given also on page~270 of \cite{shapesAndGeometry}.
Since the argumentation there contains a small gap, we
provide a proof here nevertheless for convenience.

Assume first that $x \in A$.  In this case,
$\d{A}(x) = 0$ such that $\abs{\d{A}(x) - \d{B}(x)} = \d{B}(x)$.
Since
\begin{equation*}
\dH(A, B) \ge \dHone{A}B = \sup_{y \in A} \d{B}(y) \ge \d{B}(x),
\end{equation*}
the estimate \eqref{dH as sup over d inequality} follows.
A similar argument can be applied if $x \in B$.
Thus, it remains to consider the case $x \in \R^n \setminus (A \cup B)$.
According to \lemmaref{inf/sup attained}, we can choose
$y \in B$ with $\d{B}(x) = \abs{x - y}$.
There also exists $z \in A$ such that $\d{A}(y) = \abs{z - y}$.
Then \eqref{dH as sup over d inequality} follows, since
\begin{equation*}
\d{A}(x) - \d{B}(x)
  \le \abs{x - z} - \abs{x - y}
  \le \abs{z - y} = \d{A}(y)
  \le \dHone{B}A \le \dH(A, B).
\end{equation*}

To show also \eqref{dH as sup over d}, let us assume, without
loss of generality, that $\dH(A, B) = \dHone{A}B$.
But since
\begin{equation}
\label{eq:dH as sup over d dHone}
\dHone{A}B
  = \sup_{x \in A} \d{B}(x)
  = \sup_{x \in A} \abs{\d{A}(x) - \d{B}(x)}
  \le \cnorm{\d{A} - \d{B}},
\end{equation}
the claim follows.
\end{proof}
\end{theorem}

\thref{dH as sup over d} forms the foundation for the remainder
of our paper:  It gives a representation of the Hausdorff distance
in terms of the distance functions.
Furthermore, it is also easy to actually \emph{evaluate} this representation
in practice.
In particular, if $\d{A}$ and $\d{B}$ are
given numerically on a grid, one can just consider
$\abs{\d{A}(x_i) - \d{B}(x_i)}$ for all grid points $x_i$.
The largest difference obtained in this way is guaranteed to be
at least a \emph{lower bound} for $\dH(A, B)$.
If the maximising point for \eqref{dH as sup over d} is not a grid point,
however, we can not expect to get \emph{equality} with the Hausdorff distance.
\secref{error estimates} will be devoted to a discussion
of the possible error introduced in this way.

It is sometimes convenient to use not the Hausdorff distance itself,
but the so-called \emph{complementary} Hausdorff distance
$\dH(\R^n \setminus A, \, \R^n \setminus B)$ instead.
(Particularly when dealing with open domains in applications.)
See, for instance, \cite{complementaryHausdorffConv}.
In this case, our assumption of compact sets is not fulfilled
any more, since the complements are unbounded if the sets themselves
are bounded.  However, one can verify that
\lemmaref{inf/sup attained} and \thref{dH as sup over d} are still
valid also for this situation.

\mysubsection{dH and sd}{Signed Distances}

We turn our focus now to \emph{signed} distance functions:
Since $\sd{A}$ and $\sd{B}$ are in $\C(\R^n)$ as well,
also $\cnorm{\sd{A} - \sd{B}}$ can be used as a distance
measure between $A$ and $B$.
See also Subsection~7.2.2 of \cite{shapesAndGeometry}.
This distance is, however, \emph{not} equal
to $\dH(A, B)$:

\begin{example}{dH not signed diff}
Let $0 < r < R$ be given, and define
$A = \Bc{R}0$, $B = A \setminus \B{r}0$.
This situation is depicted in \figref{dH not signed diff}.
Then $\dH(A, B) = \dHone{A}B = r$, as highlighted
in the sketch with the red line.
On the other hand, $\sd{A}(0) = -R$ while $\sd{B}(0) = r$.
Hence,
\begin{equation*}
\cnorm{\sd{A} - \sd{B}} \ge \abs{\sd{A}(0) - \sd{B}(0)} = R + r > r = \dH(A, B).
\end{equation*}

\includefig{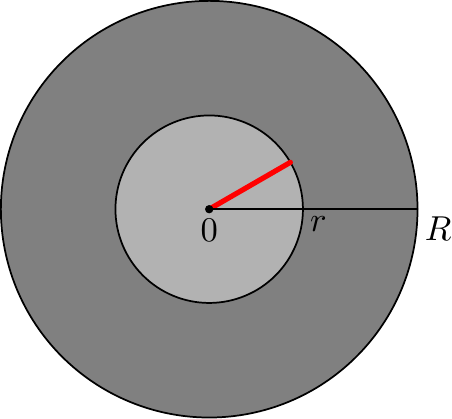}{width=0.3\textwidth}{dH not signed diff}
  {Sketch for the situation of \exref{dH not signed diff}.
   The dark ring is $B = A \cap B$, while the inner circle is
   $A \setminus B$.  The red line shows the Hausdorff distance
   between $A$ and $B$.}
\end{example}

In fact, one can show that $\cnorm{\sd{A} - \sd{B}}$ induces a \emph{stronger}
metric between the sets than either the complementary or the ordinary
Hausdorff distance alone:

\begin{theorem}{sd dist stronger}
Let $x \in \R^n$.  Then
\begin{equation}
\label{eq:sd dist stronger single x}
\abs{\sd{A}(x) - \sd{B}(x)}
 = \abs{\d{A}(x) - \d{B}(x)}
    + \abs{\d{\R^n \setminus A}(x) - \d{\R^n \setminus B}(x)}.
\end{equation}
Consequently, also
\begin{equation}
\label{eq:sd dist stronger ineq}
\max\left(
    \dH(A, B), \,
    \dH(\R^n \setminus A, \, \R^n \setminus B)
  \right)
  \le \cnorm{\sd{A} - \sd{B}}
  \le \dH(A, B) + \dH(\R^n \setminus A, \, \R^n \setminus B).
\end{equation}

\begin{proof}
Choose $x \in \R^n$ arbitrary.
If $x \in A \cap B$, then
\begin{equation*}
\abs{\d{A}(x) - \d{B}(x)} = 0, \;\;
\abs{\d{\R^n \setminus A}(x) - \d{\R^n \setminus B}(x)}
  = \abs{\sd{A}(x) - \sd{B}(x)}.
\end{equation*}
This implies the claim.
For $x \in A \setminus B$ instead, we get
\begin{equation*}
\abs{\d{A}(x) - \d{B}(x)} = \d{B}(x), \;\;
\abs{\d{\R^n \setminus A}(x) - \d{\R^n \setminus B}(x)}
  = \d{\R^n \setminus A}(x), \;\;
\abs{\sd{A}(x) - \sd{B}(x)}
  = \d{\R^n \setminus A}(x) + \d{B}(x).
\end{equation*}
Taking these together, we see that the claim is satisfied also in this case.
The two remaining cases can be handled with analogous arguments.
The relation \eqref{sd dist stronger ineq} follows by taking the
supremum over $x \in \R^n$ in \eqref{sd dist stronger single x}.
\end{proof}
\end{theorem}

Unfortunately, equality does not hold in general for
the right part of \eqref{sd dist stronger ineq}.
This is due to the fact that
taking the supremum in \eqref{sd dist stronger single x}
may yield \emph{different} maximisers for
$\abs{\d{A}(x) - \d{B}(x)}$
and $\abs{\d{\R^n \setminus A}(x) - \d{\R^n \setminus B}(x)}$.
One can also construct a simple example where this is, indeed, the case:

\begin{example}{sd dist stronger no eq}
Choose $A = [0, 1]$ and $B = [0, 3]$.
Then $\dH(A, B) = 2$, while $\dH(\R \setminus A, \, \R \setminus B) = 3/2$.
For the signed distance functions, we have $\cnorm{\sd{A} - \sd{B}} = 2$.
See also \figref{sd dist stronger no eq}, which sketches this situation.

\includesubfig{sd dist stronger no eq}
  {The situation of \exref{sd dist stronger no eq}.  The red line highlights
   the Hausdorff distance between $A$ and $B$ (left), as well as their
   complements (right).}
  {
    \subfig{0.4}{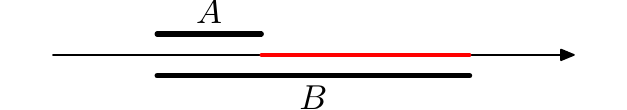}{AB}{$A$ and $B$ themselves.}
    \subfig{0.4}{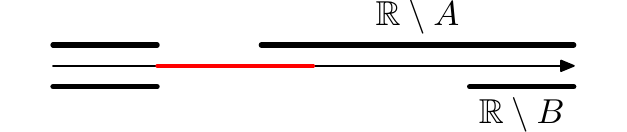}{complements}
      {$\R \setminus A$ and $\R \setminus B$.}
  }
\end{example}

That $\cnorm{\sd{A} - \sd{B}}$ is strictly stronger
than $\dH(A, B)$ also manifests itself in the
\emph{induced topology} on the space of compact subsets of $\R^n$:

\begin{example}{not same topology}
Let $A = [-1, 1]$.  For $k \in \N$, we define
\begin{equation*}
A_k = A \setminus \left(-\frac1k, \frac1k\right).
\end{equation*}
This defines a compact set $A_k \subset A$ for each $k$.
Furthermore, $\dH(A, A_k) = 1/k \to 0$ as $k \to \infty$.
In other words, $A_k \to A$ in the Hausdorff distance.
However, $\sd{A}(0) = -1$ while $\sd{A_k}(0) = 1/k$.
In particular, $\sd{A_k} \not\to \sd{A}$.
\end{example}

\exref{not same topology} implies also
that the reverse of \eqref{sd dist stronger ineq},
\begin{equation*}
\cnorm{\sd{A} - \sd{B}} \le C \cdot \dH(A, B),
\end{equation*}
can not hold for any constant $C$.
Thus, one really needs \emph{both} the ordinary \emph{and} the complementary
Hausdorff distance to get an upper bound on $\cnorm{\sd{A} - \sd{B}}$.
In other words, $\dH(A, B)$ and $\cnorm{\sd{A} - \sd{B}}$
are \emph{not} equivalent metrics.
Compare also Example~2 in \cite{complementaryHausdorffConv}:
There, it is shown that the topologies induced by the ordinary
and the complementary Hausdorff distance are not the same.
This is done with a construction similar to \exref{not same topology}.

\mysubsection{max dist}{The Maximum Distance Function}

In the final part of this section, we would like to introduce
another lower bound for $\dH(A, B)$.  This additional bound may
improve the approximation of $\dH(A, B)$ if we are not able to maximise over
all $x \in \R^n$ but only, for instance, grid points.
However, we ultimately come to the conclusion
that this bound is probably not very useful for a practical computation of
$\dH(A, B)$.
This will be discussed further at the end of the current subsection.
Hence, we will not make use of the results here in the later
\secref{error estimates}.
Since the concepts are, nevertheless,
interesting at least from a theoretical point of view,
we still give a brief presentation here.
As far as we are aware, these results have not been discussed
in the literature before.

Our initial motivation is the following:
We have seen in \thref{dH as sup over d} that the Hausdorff distance
$\dH(A, B)$ can be expressed as $\cnorm{\d{A} - \d{B}}$.
On the other hand, $\cnorm{\sd{A} - \sd{B}}$ gives \emph{not} the
Hausdorff distance.
If we are given $\sd{A}$ and $\sd{B}$ for the computation, this is unfortunate.
While it is, of course, trivial to get $\d{A}$ and $\d{B}$ from the signed
distance functions, this process throws away valuable information.
In particular, the information from the signed distance functions at points
\emph{inside} the sets can not be used.
By defining yet another type of ``distance function'', which now gives
the \emph{maximal} distance to any point in a set, we get rid
of this qualitative difference between interior and exterior points:

\begin{definition}{md}
Let $\Omega \subset \R^n$ be bounded.  The \emph{maximum distance function}
of $\Omega$ is then
\begin{equation}
\label{eq:md}
\md\Omega(x) = \sup_{y \in \Omega} \abs{x - y}.
\end{equation}
Since $\Omega$ is bounded, this is well-defined for any $x \in \R^n$.
If $\Omega$ is compact in addition, an analogous result to
\lemmaref{inf/sup attained} holds.
\end{definition}

Indeed, $\md\Omega$ is always non-negative (assuming $\Omega \ne \emptyset$).
For $\md\Omega(x)$, it does not immediately matter whether
$x \in \Omega$ or not.
Furthermore, also the maximum distance function
gives a lower bound on the Hausdorff distance,
similar to \eqref{dH as sup over d inequality}:

\begin{theorem}{md lower bound dH}
Let $A, B \subset \R^n$ be compact and choose $x \in \R^n$ arbitrarily.
Then
\begin{equation*}
\abs{\md{A}(x) - \md{B}(x)} \le \dH(A, B).
\end{equation*}

\begin{proof}
The proof is similar to the proof of \thref{dH as sup over d}:
Let $x \in \R^n$ be given.  There exist $y \in B$ with
$\md{B}(x) = \abs{x - y}$ and $z \in A$ with $\d{A}(y) = \abs{y - z}$.
Note that $\dH(A, B) \ge \abs{y - z}$ and
$\md{A}(x) \ge \abs{x - z}$.  Thus
\begin{equation*}
\md{B}(x) - \md{A}(x)
  \le \abs{x - y} - \abs{x - z} \le \abs{y - z} \le \dH(A, B).
\end{equation*}
This completes the proof if we apply the same argument
also with the roles of $A$ and $B$ exchanged.
\end{proof}
\end{theorem}

Unfortunately, the analogue of \eqref{dH as sup over d} does \emph{not} hold.
In fact, it is possible that $\md{A} = \md{B}$ everywhere on $\R^n$ but the
sets $A$ and $B$ are quite dissimilar.
Such a situation is depicted in \figref{md equal}.
Due to the ``outer ring'', which is part of $A \cap B$,
the maximum in \eqref{md} is always achieved with some $y$ from this ring.
A typical situation is shown with the point $x$ and the red line,
which highlights its maximum distance to both $A$ and $B$.
Consequently, $\md{A} = \md{B}$ and the differences between $A$ and $B$
inside the ring are not ``seen'' by the maximum distance functions at all.
Thus, we have to accept that $\cnorm{\md{A} - \md{B}} < \dH(A, B)$
can be the case.

\includefig{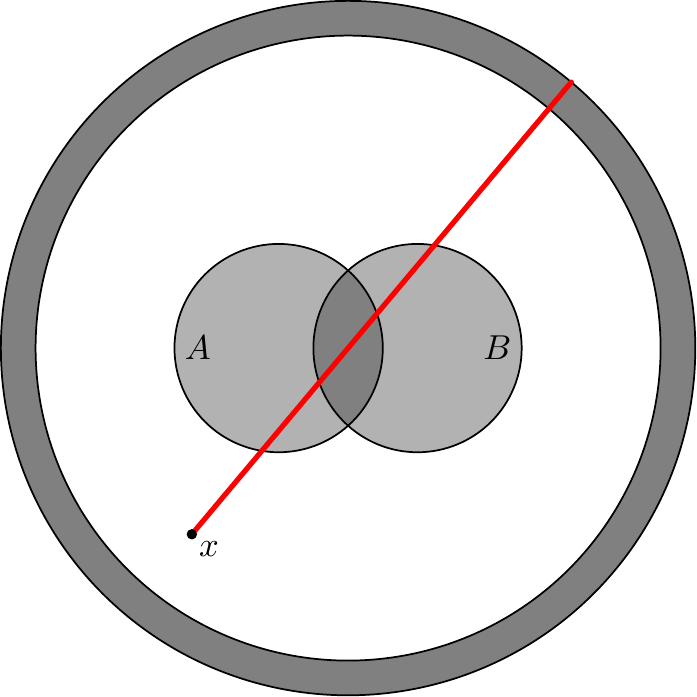}{width=0.4\textwidth}{md equal}
  {In the shown situation, $\md{A} = \md{B}$ while $A$ and $B$ are clearly
   not the same sets.  Dark regions are $A \cap B$, while lighter regions
   are $A$ or $B$ alone.
   The red line highlights the maximum distance $\md{A}(x) = \md{B}(x)$
   for some generic point $x$.}

The situations where the maximum distance functions actually carry
valuable information (as opposed to \figref{md equal}) are actually
similar to those characterised in \defref{ext dH}.
For such situations, the additional information in $\md{A}$ and $\md{B}$
could, indeed, be used to improve the approximation of $\dH(A, B)$.
However, as we will see below in \thref{ext dH estimate},
those are also the situations
where \eqref{dH as sup over d} alone already gives a very close
estimate of $\dH(A, B)$.
In these cases, we are not really in need of additional information.
On the other hand, for situations like \figref{md equal}
also the approximation of $\dH(A, B)$ from grid points is actually
difficult and extra data would be very desirable.
But particularly for those situations, the maximum distance functions
do not provide any extra data!
Furthermore, it is not clear how $\md{A}$ and $\md{B}$ can actually
be computed from, say, the level-set functions of $A$ and $B$.
It seems plausible that those functions are the viscosity solutions
of an equation similar to the Eikonal equation, and so it may be possible
to develop either a Fast Marching scheme or some other numerical method.
However, since we have just argued that we do not expect a real
benefit from the usage of the maximum distance functions in practice,
the effort involved seems not worthwhile.
For the remainder of this paper, we will thus concentrate on
\thref{dH as sup over d} as the sole basis
for our numerical computation of $\dH(A, B)$.

\mysection{error estimates}{Estimation of the Error on a Grid}

With the basic theoretical background of \secref{properties},
let us now consider the situation on a grid.
In particular, we assume that we have a rectangular, bounded grid
in $\R^n$ with uniform spacing $h$ in each dimension.
(While it is possible to generalise some of the results to non-uniform
grids in a straight-forward way, we assume a uniform spacing for simplicity.)
We denote the finite set of all grid points by $N$,
and the set of all grid cells by $C$.  For each cell $c \in C$,
$N(c)$ is the set of all grid points that span the cell (i.~e., its corners).
For example, for a $k \times k$ grid in $\R^2$ that extends
from the origin into the first quadrant, we have
\begin{equation*}
\begin{split}
N & = \set{x_{ij}}{\intrange{i,j}0{k-1}, \; x_{ij} = (i, j) h}, \\
C & = \set{c_{ij}}{\intrange{i,j}1{k-1}}, \;\;
  N(c_{ij}) = \directset{x_{i-1,j-1}, \, x_{i,j-1}, \, x_{ij}, \, x_{i-1,j}}.
\end{split}
\end{equation*}
Let us assume that we know the distance functions of $A$ and $B$
on each grid point, i.~e., $\d{A}(x)$ and $\d{B}(x)$ for all $x \in N$.
We furthermore assume that these values are known without approximation error.
This is, of course, not realistic in practice.  However, the approximation error
in describing the geometries and computing their distance functions is
a matter outside the scope of this paper.
Finally, let us also assume that the grid is large enough to cover
the sets.  In particular:
For each $y \in A \cup B$, there should exist a grid cell $c \in C$
such that $y$ is contained in the convex hull $\convex{N(c)}$ of
the corners of $c$.
If this is not the case, the grid is simply inadequate to capture
the geometrical situation.

\mysubsection{worst case}{Worst-Case Estimates}

In order to approximate $\dH(A, B)$ from the distance functions on
our grid, we make use of \eqref{dH as sup over d}.
In particular, we propose the following straight-forward approximation:
\begin{equation}
\label{eq:dH approx}
\dH(A, B) \approx \dHapprox(A, B) = \max_{x \in N} \abs{\d{A}(x) - \d{B}(x)}
\end{equation}
From \eqref{dH as sup over d inequality}, we know that this is, at least,
an exact lower bound.  However, in the general case, an approximation error
\begin{equation*}
0 \le \delta = \abs{\dH(A, B) - \dHapprox(A, B)} = \dH(A, B) - \dHapprox(A, B)
\end{equation*}
will be introduced by using \eqref{dH approx}.
This is due to the fact that we only maximise over grid points.
The real maximiser of the supremum in \eqref{dH as sup over d},
on the other hand, may not be a grid point.


Let us now analyse the approximation error $\delta$.
We have seen in the proof of \thref{dH as sup over d} that
\begin{equation*}
\dH(A, B) = \sup_{y \in A \cup B} \abs{\d{A}(y) - \d{B}(y)}
  = \max_{c \in C} \sup_{y \in \convex{N(c)}} \abs{\d{A}(y) - \d{B}(y)}.
\end{equation*}
Note that this is still an exact representation, with no approximation
error introduced so far.  We have just split up the supremum
over $A \cup B$ into grid cells, but we still take into account
\emph{all} points contained in a grid cell, not just its corners.
This is achieved by using the convex hull $\convex{N(c)}$ instead of
the finite set $N(c)$ alone.
On the other hand, the approximation \eqref{dH approx} can be formulated as
\begin{equation*}
\dHapprox(A, B) = \max_{c \in C} \max_{x \in N(c)} \abs{\d{A}(x) - \d{B}(x)}.
\end{equation*}
Comparing both equations, we see that the approximation error $\delta$
is introduced precisely by the step from $\convex{N(c)}$ to $N(c)$.
We can now formulate and prove a very general \emph{upper bound}
on $\dH(A, B)$:

\begin{theorem}{general upper bound}
Let $x \in N$ be a grid point and $y \in \R^n$ be arbitrary.  We set
\begin{equation*}
t(x, y) = \begin{alternatives}
            \alternative{\abs{x - y}}{x \in A,}
            \alternative{2 \abs{x - y}}{x \not\in A.}
          \end{alternatives}
\end{equation*}
For a cell $c \in C$, we define furthermore
\begin{equation}
\label{eq:general upper bound cell}
\overline{d}(c)
  = \sup_{y \in \convex{N(c)}} \min_{x \in N(c)}
      \left(\abs{\d{A}(x) - \d{B}(x)} + t(x, y)\right).
\end{equation}
Then
\begin{equation*}
\dHone{A}B \le \max_{c \in C'(A)} \overline{d}(c).
\end{equation*}
Here, $C'(A) = \set{c \in C}{\convex{N(c)} \cap A \ne \emptyset}$
is the set of all grid cells which contain some part of $A$.

Similarly, $\dHone{B}A$ and thus $\dH(A, B)$ can be estimated.

\begin{proof}
We will show that
\begin{equation*}
\sup_{y \in A \cap \convex{N(c)}} \abs{\d{A}(y) - \d{B}(y)}
  \le \overline{d}(c)
\end{equation*}
for each $c \in C$.
The claim then follows from \eqref{dH as sup over d dHone}.
So choose $c \in C$, $y \in A \cap \convex{N(c)}$ and $x \in N(c)$.
It remains to verify that
\begin{equation*}
\abs{\d{A}(y) - \d{B}(y)} = \d{B}(y) \le \abs{\d{A}(x) - \d{B}(x)} + t(x, y).
\end{equation*}
Assume first that $x \in A$.
Since $\d{B}$ has Lipschitz constant one, we really get
\begin{equation*}
\abs{\d{A}(x) - \d{B}(x)} + t(x, y)
  = \d{B}(x) + \abs{x - y} \ge \d{B}(y)
\end{equation*}
in this case.
Assume now $x \not\in A$.  Since $y \in A$ and thus
$\d{A}(y) = 0$, Lipschitz continuity of $\d{A}$ implies that
$\d{A}(x) \le \abs{x - y}$.
Using this auxiliary result, we get that also in this case
\begin{equation*}
\abs{\d{A}(x) - \d{B}(x)} + t(x, y)
  \ge \d{B}(x) - \d{A}(x) + 2 \abs{x - y}
  \ge \d{B}(x) + \abs{x - y} \ge \d{B}(y).
\end{equation*}
Hence, the claim is shown.
\end{proof}
\end{theorem}

Even though the formulation of \thref{general upper bound} is complicated,
the idea behind it is quite simple:
Since the distance functions are Lipschitz continuous, also
the function $\abs{\d{A} - \d{B}}$, which we have to maximise over
$\convex{N(c)}$ for each grid cell, is Lipschitz continuous.
This allows us to estimate the maximum in terms of the function's
values at the corners (which are known).
We are even allowed to try all corners and use
the smallest resulting upper bound.  This is what happens in
\eqref{general upper bound cell}.
Furthermore, the Lipschitz constant depends on whether or not the corner
is in $A$.  (If it is, $\d{A}$ vanishes, which reduces
the Lipschitz constant to just that of $\d{B}$.  Otherwise, we have to
use two as the full Lipschitz constant of $\abs{\d{A} - \d{B}}$.)
This is the role that $t(x, y)$ plays.
It gives the ``distance'' between $x$ and $y$ based on the applicable
Lipschitz constant.


Coupled with the fact that $\dHapprox(A, B)$ is a lower bound
for the exact Hausdorff distance, the upper bound in \thref{general upper bound}
allows us now to estimate $\delta$.
However, evaluating \eqref{general upper bound cell} is
difficult and expensive in practice (although it can be done in theory).
Hence, we will now draw some conclusions that simplify the
upper bound.  As a first result, let us consider the worst case
where $x \not\in A$ for all corners $x \in N(c)$ of some cell:

\begin{corollary}{worst case bound}
\thref{general upper bound} implies for the error estimate:
\begin{equation*}
\delta \le \sqrt{n} \cdot h
\end{equation*}

\begin{proof}
Let $c \in C$ be some grid cell and $y \in \convex{N(c)}$.
Then there exists $x \in N(c)$ such that
\begin{equation*}
t(x, y) \le 2 \abs{x - y} \le 2 \cdot \frac{\sqrt{n} \cdot h}2
  = \sqrt{n} \cdot h.
\end{equation*}
This is simply due to the fact that the grid cell's longest diagonal
has length $\sqrt{n} \cdot h$.  Consequently, in the worst case the nearest
corner $x$ has half that distance to $y$.
Hence also
\begin{equation*}
\begin{split}
\min_{x \in N(c)} \left(\abs{\d{A}(x) - \d{B}(x)} + t(x, y)\right)
  & \le \max_{x \in N(c)} \abs{\d{A}(x) - \d{B}(x)}
          + \min_{x \in N(c)} t(x, y) \\
  & \le \max_{x \in N(c)} \abs{\d{A}(x) - \d{B}(x)} + \sqrt{n} \cdot h.
\end{split}
\end{equation*}
This estimate can be used for
$\overline{d}(c)$ from \eqref{general upper bound cell}.
Consequently, \thref{general upper bound} implies
\begin{equation*}
\dHone{A}B
  \le \max_{c \in C} \max_{x \in N(c)} \abs{\d{A}(x) - \d{B}(x)}
        + \sqrt{n} \cdot h
   =  \dHapprox(A, B) + \sqrt{n} \cdot h.
\end{equation*}
Since the same estimate also works for $\dHone{B}A$, the claim follows.
\end{proof}
\end{corollary}


Taking a closer look, though, the worst-case situation considered
above is quite strange.  In principle, it \emph{can} happen that there is some
cell $c$ with $A \cap \convex{N(c)} \ne \emptyset$ but for which
all corners are not in $A$.
Such a situation is depicted in \figref{all corners not in A}.
However, in practice such a case is very unlikely to occur.
In particular, assume that we describe the set $A$ by a level-set function
$\phi$, and that $\phi(x) > 0$ for all corners $x \in N(c)$
of some grid cell $c$.
In that case, there is no way of knowing whether, in reality, there is some part
of $A$ inside the cell or not.
The grid is simply too coarse to ``see'' such geometric details.
Consequently, it makes sense to assume the simplest possible situation,
namely that $A \cap \convex{N(c)} = \emptyset$ for all such cells $c$.
Thus, we make the following additional assumption:

\begin{definition}{suitable grid}
Consider grid cells $c \in C$ such that $x \not\in A$
for all corners $x \in N(c)$.
If $\convex{N(c)} \cap A = \emptyset$ for all those $c$,
the grid is said to be \emph{suitable} for $A$.
\end{definition}

\includefig{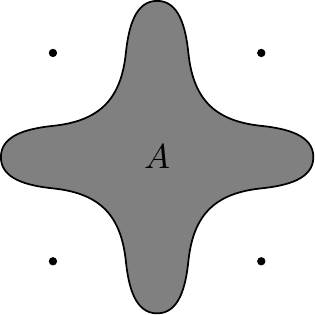}{width=0.2\textwidth}{all corners not in A}
  {A possible pathological situation where the grid is not suitable for $A$.}

In the case of a suitable grid (for both $A$ and $B$), we get
the ``reduced Lipschitz constant'' in \eqref{general upper bound cell}
for at least one corner per relevant cell.
This allows us to lower the error estimate:

\begin{corollary}{bound for suitable grid}
Let the grid be suitable for $A$ and $B$.
Furthermore, we introduce the dimensional constant
\begin{equation}
\label{eq:bound dim constant}
\Delta_n = \sup_{y \in Q}
            \min\left(\abs{y}, \, \min_{x \in N'} 2 \abs{x - y}\right).
\end{equation}
Here, $Q = \set{y \in \R^n}{0 \le y_i \le 1 \text{ for all $\intrange{i}1n$}}$
is the unit square, and
\begin{equation*}
N' = \set{x \in \R^n}
         {x_i \in \directset{0, 1} \text{ for all $\intrange{i}1n$}}
      \setminus \directset{0}
\end{equation*}
is the set of its corners except for the origin.  Then,
\begin{equation*}
\delta \le \Delta_n \cdot h.
\end{equation*}

\begin{proof}
Let $c \in C$ be a cell with $\convex{N(c)} \cap A \ne \emptyset$.
Since the grid is assumed to be suitable, we know that there exists
at least one corner $x_0 \in N(c)$ with $x_0 \in A$.
Hence, for arbitrary $y \in \convex{N(c)}$,
\begin{equation*}
\min_{x \in N(c)} t(x, y)
  \le \min\left(\abs{y - x_0}, \,
                \min_{x \in N(c) \setminus \directset{x_0}}
                  2\abs{y - x}\right)
  \le \Delta_n \cdot h.
\end{equation*}
With this, the claim follows as in the proof of \corref{worst case bound}.
\end{proof}
\end{corollary}


The most difficult part of \corref{bound for suitable grid}
is probably the strange dimensional constant $\Delta_n$ defined
in \eqref{bound dim constant}.
This constant replaces the functions $t(x, \cdot)$ for $x \in N(c)$.
It can be interpreted like this:
Let spherical fronts propagate starting from all corners of the
unit square $Q$.  The front starting at the origin has speed one,
while the other fronts have speed $1/2$.
Over time, the fronts will hit each other, and will reach all parts of $Q$.
The value of $\Delta_n$ is precisely the
time it takes until all points in $Q$ have been hit by at least one front.
For the case $n = 2$, these arrival times are shown in \subfigref{delta}2.
The correct value of $\Delta_2$ is the maximum attained at both
spots with the darkest red (one in the north and one in the east).
\subfigref{delta}3 shows the maximising points (red and black)
over the unit cube for $n = 3$.  Since the expression that is maximised
in \eqref{bound dim constant} is symmetric with respect to permutation of
the coordinates, there are six maximisers.  The highlighted one
sits at the intersection of the spheres originating from the three
corners marked in blue.
Based on these observations and some purely geometrical arguments,
one can calculate
\begin{equation*}
\Delta_1 = \frac23 \approx \num{0.67}, \;\;
\Delta_2 = \frac23 \sqrt{5 - \sqrt7} \approx \num{1.02}, \;\;
\Delta_3 = \frac23 \sqrt{8 - \sqrt{19}} \approx \num{1.27}.
\end{equation*}
These constants are a clear improvement over the estimate of
\corref{worst case bound}.
In fact, the bound in \corref{bound for suitable grid} is actually sharp.
To demonstrate this, we will conclude this subsection with
an example in two dimensions that really attains the maximal
error $\delta = \Delta_2 \cdot h$:

\includesubfig{delta}
  {The computation of $\Delta_2$ and $\Delta_3$ from the maximisation
   of the arrival times over $Q$.}
  {
    \subfig{0.8}{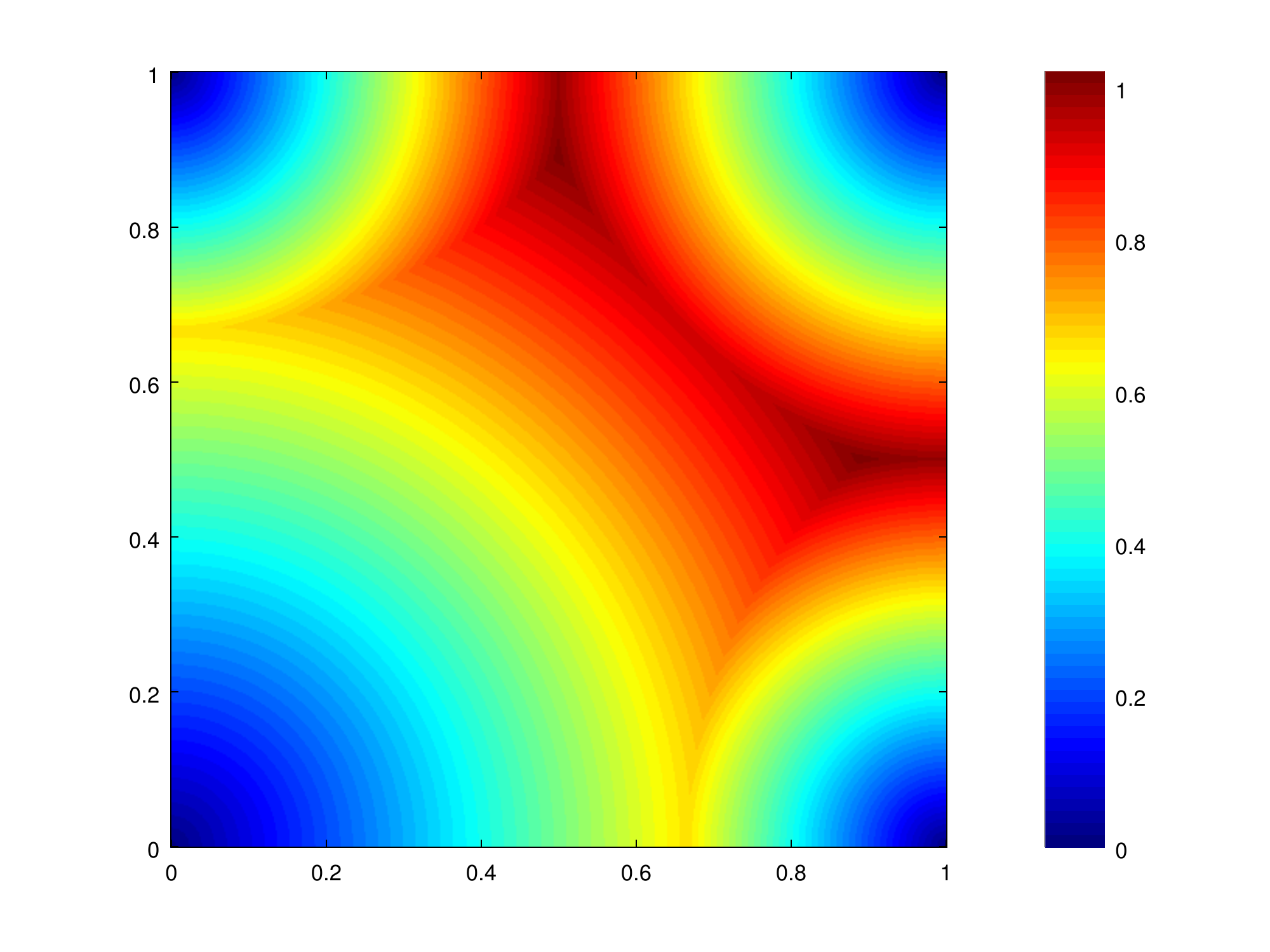}2{Arrival times for $\Delta_2$.} \\
    \subfig{0.8}{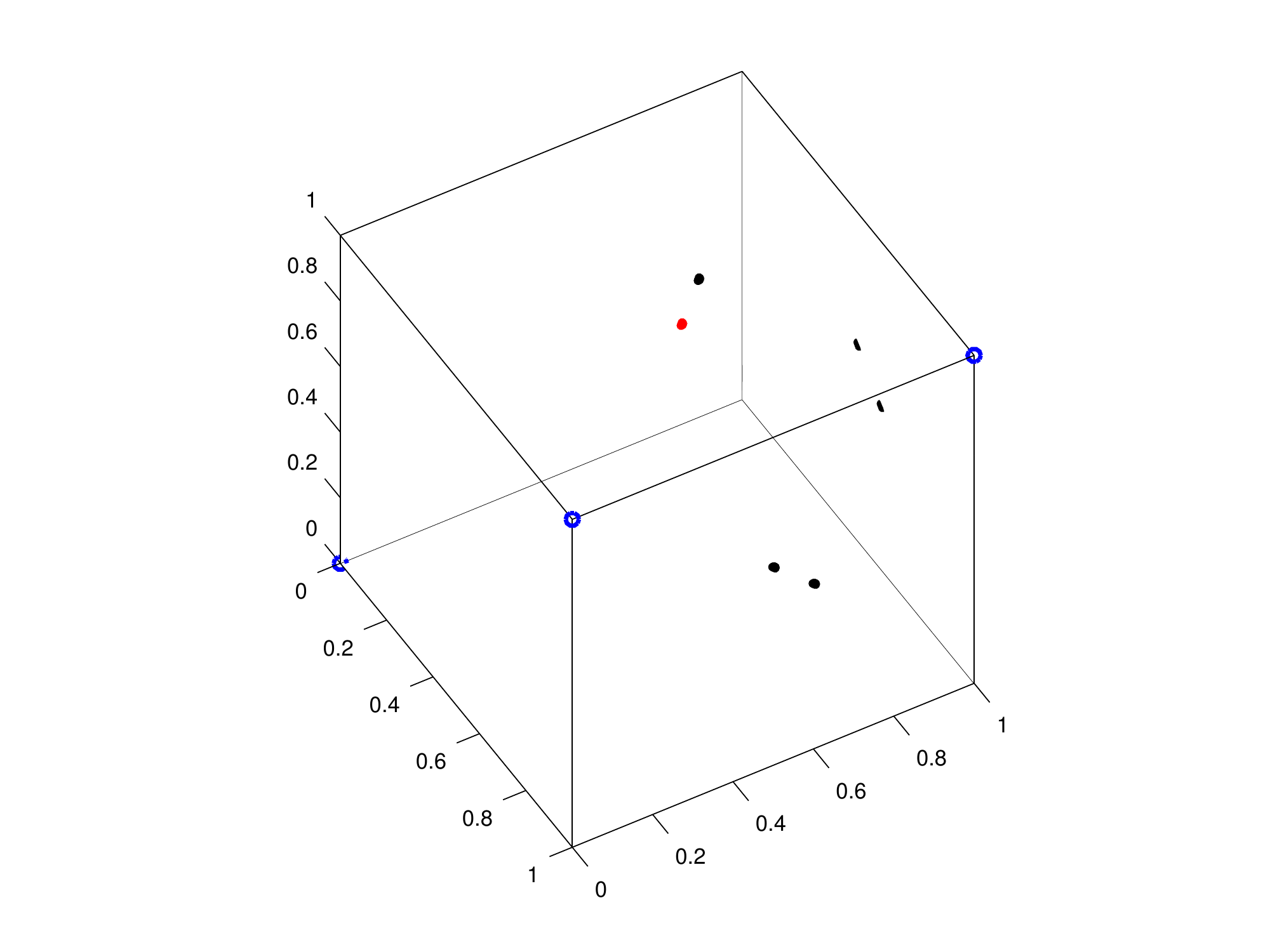}3{Maximising points for $\Delta_3$.} \\
  }


\begin{example}{maximal error}
For simplicity, assume $h = 1$.
We consider the situation sketched in \figref{maximal error}.
Observe first that all grid points except $a$, $b$, $c$ and $d$
are part of $A \cap B$, and thus $\d{A} = \d{B} = 0$ for them.
Consequently, we only have to consider these four points
in order to find $\dHapprox(A, B)$.  For symmetry reasons,
it is actually enough to concentrate only on $a$ and $b$.
The point $p$ corresponds to the position with maximal arrival time,
as seen also in \subfigref{delta}2.  It is characterised
by requiring
\begin{equation}
\label{eq:maximal error p}
\abs{a - p} = \abs{d - p} = \frac{r}2, \;\;
\abs{b - p} = \abs{c - p} = r.
\end{equation}
Solving these equations for the coordinates of $p$ and the radius $r$ yields
\begin{equation*}
p = \left(\frac{8 - \sqrt7}6, \; \frac12\right), \;\;
r = \frac23 \sqrt{5 - \sqrt7}.
\end{equation*}
Note specifically that $r = \Delta_2$.
(In fact, a construction similar to this one can be used to calculate
$\Delta_2$ in the first place.)
The relations \eqref{maximal error p} can also be seen in the sketch:
The dotted circle has radius $r$ and centre $p$.
The points $b$ and $c$ lie on it.
The two smaller circles (which define the exclusion from $A$)
have centres $a$ and $d$ with radius $r/2$, and $p$ lies on both
of them.

$B$ is chosen in such a way that $p$ is also at the centre
of its hole.  Consequently, $p$ is the point in $A$ that achieves
$\dH(A, B) = \dHone{A}B = \d{B}(p)$.  This is indicated by the red line.
Let us also introduce $\rho$ as the width of the small ring
between the dotted circle and the dark region.
Then $\rho = \d{B}(b)$ and $\dH(A, B) = r + \rho$.
Furthermore, note that
\begin{equation*}
\d{A}(a) = \frac{r}2, \;\; \d{B}(a) = \frac{r}2 + \rho.
\end{equation*}
Hence,
\begin{equation*}
\dHapprox(A, B) = \abs{\d{A}(a) - \d{B}(a)} = \abs{\d{A}(b) - \d{B}(b)}
  = \rho.
\end{equation*}
This also implies that $\delta = r = \Delta_2 \cdot h$, which is, indeed, the
largest possible bound permitted by \corref{bound for suitable grid}.

\includefig{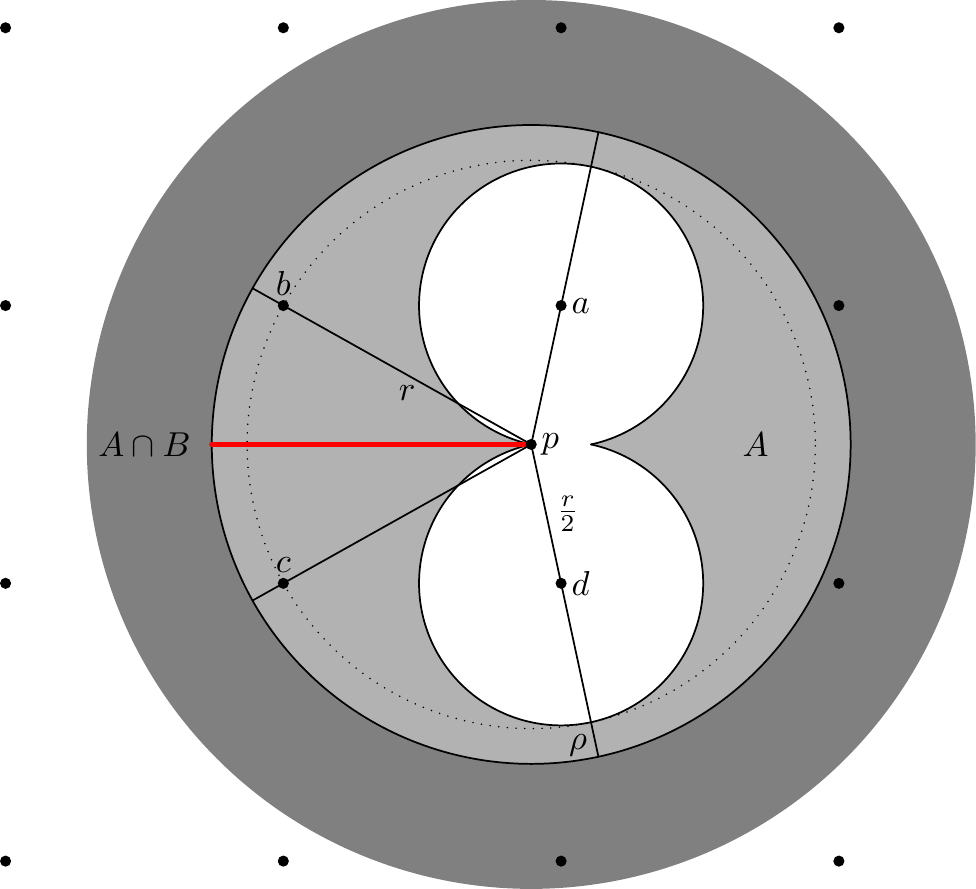}{width=0.7\textwidth}{maximal error}
  {The situation considered in \exref{maximal error}.  The dark region
   is $A \cap B$, the lighter part inside the circle is $A \setminus B$.
   The red line indicates the Hausdorff distance $\dH(A, B) = \dHone{A}B$.}
\end{example}

\mysubsection{good case}{External Hausdorff Distances}

As we have promised, the situation from \exref{maximal error}
shows that one can not, in general, expect a better error
estimate than \corref{bound for suitable grid}.
However, considering \figref{maximal error}, we also observe
that the situation there is quite strange.
Thus, there is hope that we can get stronger estimates if we
add some more assumptions on the geometrical situation.
This is the goal of the current subsection.  It will turn out
in \thref{ext dH estimate} that this is, indeed, possible.
Consider, for example, \subfigref{ext dH}{grid}:
There, $\dH(A, B) = \dHone{A}B = \abs{x - y}$.
Furthermore, all points $z$ along the external black line
satisfy $\abs{\d{A}(z) - \d{B}(z)} = \abs{x - y}$.
Consequently, \emph{all those points} are maximisers of
\eqref{dH as sup over d}.
If a grid point happens to lie somewhere on this line,
$\dHapprox(A, B)$ is exact.  But even if this is not the case
(as shown in the figure), $\abs{\d{A} - \d{B}}$ will be very close
to $\abs{x - y}$ for grid points that are far away from the sets and
close to the line.  In all of these cases, we can expect
$\dHapprox(A, B)$ to be much closer to $\dH(A, B)$ than the bounds
from the previous \subsecref{worst case} tell us.
Furthermore, the estimate will be more precise the further away we can
go on the external line.
Two conditions determine how far that really is:
First, of course, the size of our finite grid is a clear restriction.
Second, we need that $\d{A}(z) = \abs{z - x}$ and $\d{B}(z) = \abs{z - y}$
for the points $z$ on the external line that we consider.
This means that $x$ and $y$ must be the closest points to $z$
of the sets $A$ and $B$, respectively.
The latter is a purely geometrical condition on $A$ and $B$, and is not
related to the grid.  Let us formalise it:

\includesubfig{ext dH}
  {The situation of an \emph{external Hausdorff distance} from \defref{ext dH}.
   The red lines indicate the Hausdorff distances
   $\dH(A, B) = \dHone{A}B$.  The dark regions are, as always, $A \cap B$,
   while the lighter parts are $A$ or $B$ alone.}
  {
    \subfig{0.45}{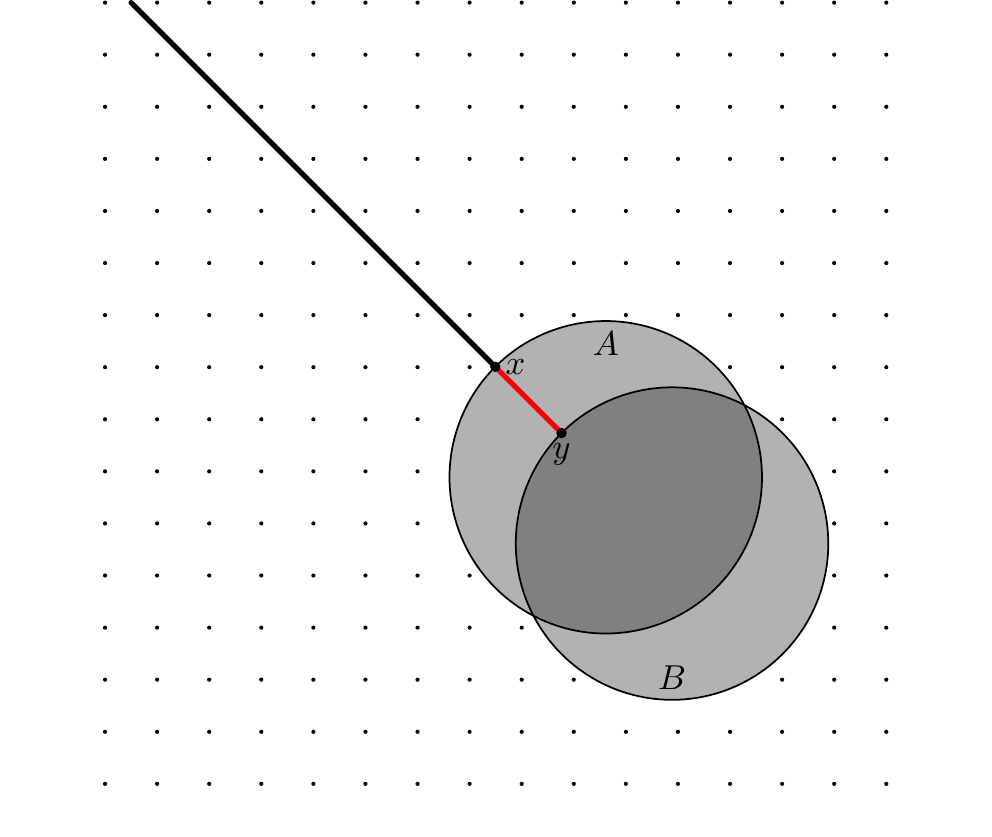}{grid}
      {The external line with $\dH(A, B) = \abs{\d{A} - \d{B}}$.}
    \subfig{0.45}{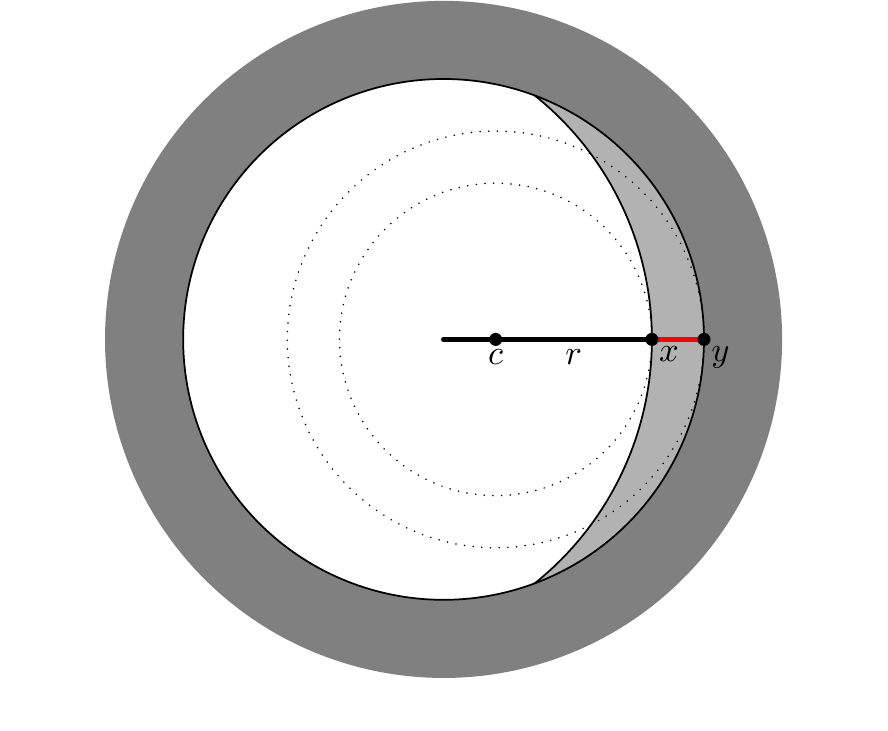}{basic}
      {External Hausdorff distance with a restricted $r$.}
  }

\begin{definition}{ext dH}
Assume that $\dH(A, B) = \dHone{A}B = \abs{x - y} > 0$
with $x \in A$ and $y \in B$.
Let $r > 0$ and set $d = (x - y) / \abs{x - y}$ as well
as $c = x + r d$ and $R = r + \dH(A, B)$.
We say that $A$ and $B$ admit an \emph{external Hausdorff distance}
with radius $r$ if
\begin{equation}
\label{eq:ext dH}
\Bc{r}c \cap A = \directset{x}, \;\;
\Bc{R}c \cap B = \directset{y}.
\end{equation}
\end{definition}

The condition in \defref{ext dH} is quite technical, but it is relatively
easy to understand and verify for concrete situations (as long as
it is known where the Hausdorff distance is attained).
It is related to the \emph{skeleton} of the sets $A$ and $B$,
for which we refer to Section~3.3 of \cite{shapesAndGeometry}.
We will see later in \corref{convex sets} that, for instance,
convex sets admit an external Hausdorff distance for arbitrary radius $r$,
and that \defref{ext dH} applies in a lot of additional practical situations.
Even for non-convex sets, an external Hausdorff
distance with some restriction on the possible $r$ may be admissible.
See, for instance, the situation in \subfigref{ext dH}{basic}.
A possible choice for $r$ and $c$ is shown there.
(The furthest possible $c$ is at the end of the black line.)
The dotted circles are $\Bc{r}c$ and $\Bc{R}c$.
One can see that the inner one only touches $A$ at $x$, and the outer one does
the same with $B$ at $y$.
This is the geometrical meaning of \eqref{ext dH}.
Due to this property, we know that $\d{A}(c) = \abs{c - x}$
and $\d{B}(c) = \abs{c - y}$.
One can also verify that the condition \eqref{ext dH}
gets strictly stronger if we increase $r$.
In other words, if an external Hausdorff distance with $r$ is admissible,
this is automatically also the case for all radii $s < r$.


Based on this concept of external Hausdorff distances,
we can now formalise the motivating argument
about better error bounds for this situation:

\begin{theorem}{ext dH estimate}
Let $A$ and $B$ admit an external Hausdorff distance
with $r > 0$.  Let $h$ be the grid spacing, and assume that
the grid is chosen large enough.
Then, for $h \to 0$,
\begin{equation}
\label{eq:ext dH estimate}
\delta \le \frac{n}2 \frac{h^2}r + \bigO{h^3}.
\end{equation}

\begin{proof}
We use the same notation as in \defref{ext dH}.
In particular, let $\dH(A, B) = \dHone{A}B = \abs{x - y}$
with $x \in A$ and $y \in B$.
We also use $c$ and $R$ as in the definition.
If $z$ is a point next to the straight line $x$--$c$,
we can project it onto this line.  Let the resulting point be called $c'$,
then $c$--$c'$--$z$ and $x$--$c'$--$z$ are right triangles.
This situation is shown in \figref{ext dH estimate}.
According to the sketch, we set
\begin{equation*}
\rho = R - \abs{z - c} = R - \sqrt{a^2 + b^2}.
\end{equation*}
Note that the dotted circle $\B\rho{z}$ is entirely contained in $\B{R}{c}$.
By \eqref{ext dH}, this implies that $\d{B}(z) \ge \rho$.
Since $x \in A$, we also know
$\d{A}(z) \le \abs{z - x} = \sqrt{b^2 + (r - a)^2}$.
Both inequalities together yield
\begin{equation*}
\dHapprox(A, B) \ge \abs{\d{A}(z) - \d{B}(z)}
  \ge \d{B}(z) - \d{A}(z) \ge R - \sqrt{a^2 + b^2} - \sqrt{b^2 + (r - a)^2}.
\end{equation*}
(Assuming that $z$ is a grid point.)
On the other hand, since we have an external Hausdorff distance, also
\begin{equation*}
\dH(A, B) = \abs{\d{A}(c') - \d{B}(c')} = \d{B}(c') - \d{A}(c')
  = (R - a) - (r - a)
\end{equation*}
holds.
Hence,
\begin{equation}
\label{eq:ext dH estimate delta}
\delta = \dH(A, B) - \dHapprox(A, B)
  \le \sqrt{b^2 + a^2} - a + \sqrt{b^2 + (r - a)^2} - (r - a).
\end{equation}

So far, $z$ was just an (almost) arbitrary grid point.
We will now try to choose it in a way that reduces the bound
on $\delta$ as much as possible.
For this, observe that \eqref{ext dH estimate delta} contains
two terms of the form $s \mapsto (\sqrt{b^2 + s^2} - s)$ and that this
function is decreasing in $s$.
Thus, in order to get a small bound, we would like to choose
both values of $s$, namely $a$ and $r - a$, as large as possible.
Consequently, we want $a \approx r/2$.
Let $m = (c + x) / 2$ be the precise midpoint between $c$ and $x$.
Since $\sqrt{n} \cdot h$ is the longest diagonal of the grid cells,
there exists a grid point $z \in N$ with $\abs{m - z} \le \sqrt{n}/2 \cdot h$.
Choosing $c'$ as the projection of $z$ onto the line $x$--$c$ as before,
this implies that
\begin{equation*}
a \ge \frac{r}2 - \abs{m - c'} \ge \frac{r - \sqrt{n} \cdot h}2, \;\;
r - a \ge \frac{r - \sqrt{n} \cdot h}2, \;\;
b = \abs{z - c'} \le \frac{\sqrt{n}}2 h.
\end{equation*}
(Since $\abs{m - z}^2 = \abs{m - c'}^2 + \abs{z - c'}^2$, not all of
these estimates can be sharp at the same time.
It may be possible to refine them and get smaller bounds below, but
we do not attempt to do that for simplicity.)
Substituting in \eqref{ext dH estimate delta}
yields
\begin{equation*}
\delta
  \le \sqrt{n h^2 + (r - \sqrt{n} \cdot h)^2} - (r - \sqrt{n} \cdot h).
\end{equation*}
Series expansion of this result for $h \to 0$ finally
implies the claimed estimate \eqref{ext dH estimate}.

\includefig{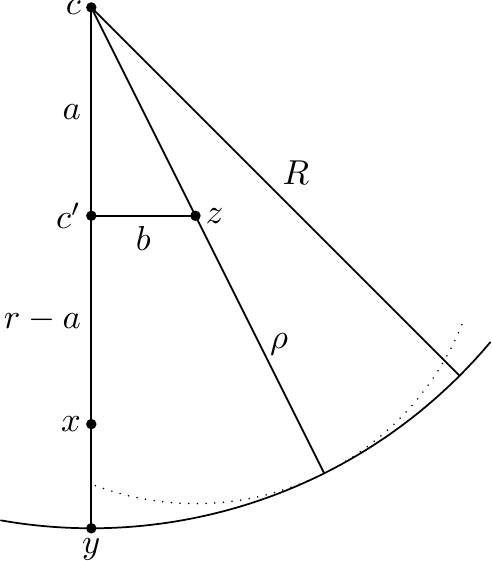}{width=0.3\textwidth}{ext dH estimate}
  {Deriving the bound \eqref{ext dH estimate delta} in the
   proof of \thref{ext dH estimate}.}
\end{proof}
\end{theorem}


A particular situation in which \defref{ext dH} is satisfied
is that of \emph{convex} sets (see \subfigref{ext dH}{grid}).
For them, \eqref{ext dH estimate}
holds with arbitrary $r$ as long as the grid is large enough
to accommodate for the far-away points:

\begin{corollary}{convex sets}
Let $A$ and $B$ be compact and convex.
Then $A$ and $B$ admit an external Hausdorff distance
for arbitrary $r > 0$.
Consequently, \eqref{ext dH estimate} applies for all $r$ for which
the grid is large enough.

\begin{proof}
We exclude the trivial case $A = B$, since \eqref{ext dH estimate}
is obviously fulfilled for that situation anyway.
Assume, without loss of generality, that $\dH(A, B) = \dHone{A}B$.
Let $x \in A$ and $y \in B$ be given with $\dHone{A}B = \abs{x - y} > 0$
according to \lemmaref{inf/sup attained}.
Choose $r > 0$ arbitrarily and let $c$ be as in \defref{ext dH}.

The assumption $\dHone{A}B = \d{B}(x) = \abs{x - y}$ means
that $y$ is the closest point in $B$ to $x$.
In other words, $\B{s}x \cap B = \emptyset$, where we have set $s = \abs{x - y}$
for simplicity.
This is depicted with the dotted circle (which is outside of $B$)
in \figref{convex sets}.
The dotted lines indicate half-planes perpendicular
to the line $x$--$y$ and through $x$ and $y$, respectively.
Assume for a moment that we have some point $z \in B$ that is
``above'' the ``lower'' half-plane.
Due to convexity of $B$, this would imply that the whole
line $z$--$y$ must be inside of $B$.  This, however,
contradicts $\B{s}x \cap B = \emptyset$ as indicated by the
red part of $z$--$y$.
Hence, the half-plane through $y$ separates $B$ and $x$.
Similarly, we can show that the half-plane through $x$ separates $A$ and $c$:
Assume that $z' \in A$ is ``above'' this half-plane.
Then $\d{B}(z') > s$ must be the case, as shown by the blue line.
But this is a contradiction, since $\d{B}(z') \le \dHone{A}B = s$
for all $z' \in A$.
These separation properties of the half-planes, however, in turn imply
\eqref{ext dH}.
Thus, everything is shown.

\includefig{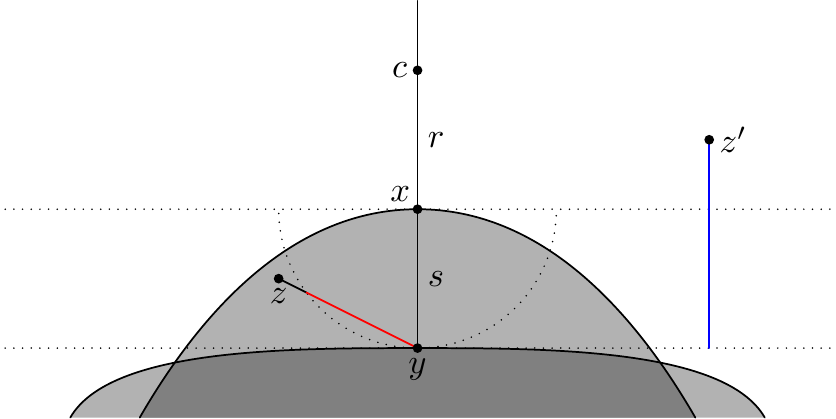}{width=0.6\textwidth}{convex sets}
  {Sketch for the proof of \corref{convex sets}.}
\end{proof}
\end{corollary}

Let us also remark that the proof of \corref{convex sets} stays
valid as long as the sets $A$ and $B$ are ``locally convex''
in a neighbourhood of $x$ and $y$.  This is an important
situation for a lot of potential applications:
We already mentioned above that our own motivation
for computing the Hausdorff distance is to measure convergence
during shape optimisation.  In this case, it is often the case that
the Hausdorff distance is already quite small in relation to the sets.
For a lot of these situations, the largest difference between the sets
is attained in a way similar to \figref{convex sets},
even if the sets themselves need not necessarily be convex.

\mysubsection{numerics}{Numerical Demonstration}

\includesubfig{circle in ring}
  {The example situation (schematically) used for \subsecref{numerics}.
   We have a dark ``outer ring'' $A \cap B$ and a smaller circle
   $A \setminus B$ that is inside.  The Hausdorff distance
   $\dH(A, B) = \dHone{A}B$ is indicated with the red line.
   In right plot, the external line with $\abs{\d{A} - \d{B}} = \dH(A, B)$
   is also indicated.}
  {
    \subfig{0.45}{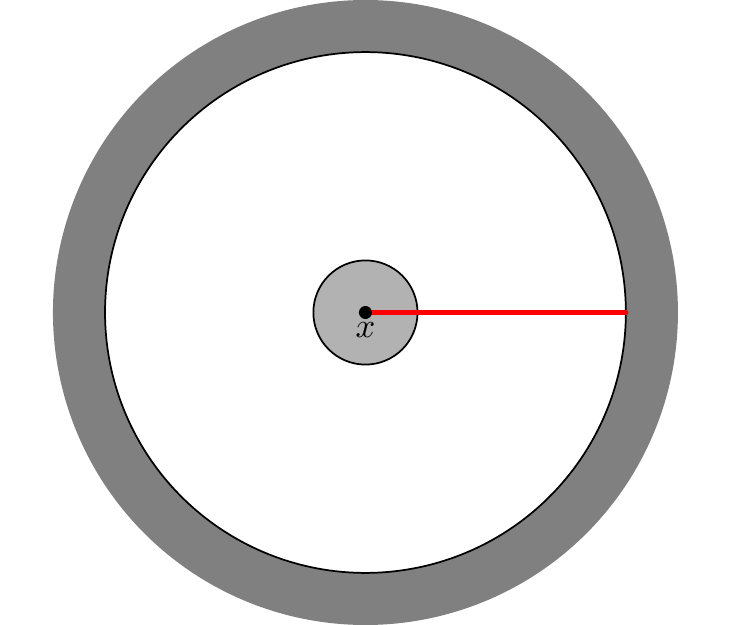}{centre}{Inner circle at the origin.}
    \subfig{0.45}{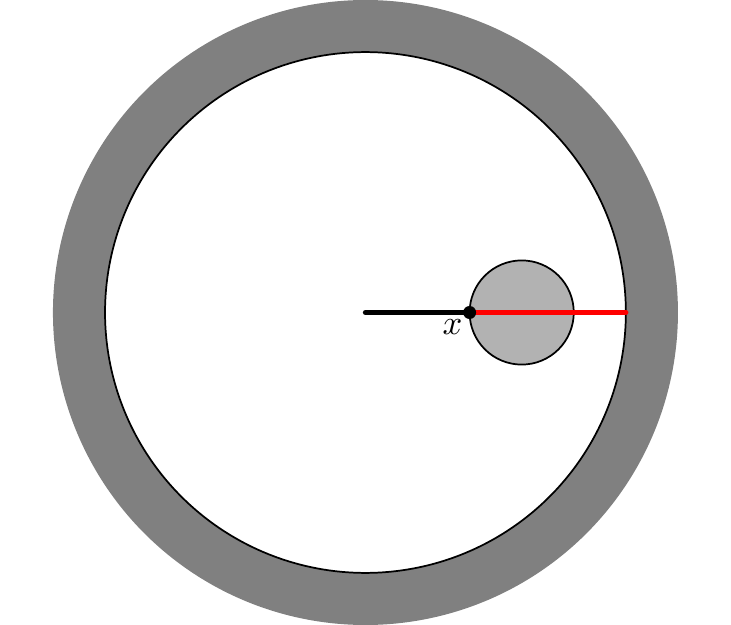}{moved}
      {Displaced inner circle.}
  }

To conclude this section,
let us give a numerical demonstration of the results presented so far.
The situation that we consider is depicted schematically
in \figref{circle in ring}:  We have an ``outer ring'' which is part
of both $A$ and $B$, and an inner circle (corresponding to $A \setminus B$)
is placed within.  Note that this is already a situation
where we have non-convex sets.  For the inner circle at the origin
as in \subfigref{circle in ring}{centre}, no external Hausdorff
distance is admissible.  This is due to the fact that the point
$x \in A$ that achieves $\dH(A, B) = \dHone{A}B = \d{B}(x)$
is in the \emph{interior} of $A$.
If we displace the inner circle, the point $x$
will be on the boundary as soon as the origin is no longer part of $A$.
In these situations, we have an external Hausdorff distance with a restricted
maximal radius $r$.  This is indicated in \subfigref{circle in ring}{moved}.
In our calculations, the outer circle has a radius of nine and the
inner circle's radius is one.
\figref{circle in ring} shows other proportions since this
makes the figure clearer.
However, \emph{qualitatively}, the situations shown are exactly those
that will be used in the following.

Let us first fix the grid spacing $h$ and consider the
effect of moving the inner circle.
The approximation error $\delta$ of the exact Hausdorff distance
is shown (in units of $h$) in \figref{plot movement}.
The blue curve shows $\delta$ under the assumption that $\sd{A}$ and $\sd{B}$
are known on the grid points without any approximation error.
This is the situation we have discussed theoretically above.
The red curve shows the error if we also compute the signed distance
functions themselves using the Fast Marching code in \cite{package}.
This is a situation that is more typical in practice,
where often only \emph{some} level-set functions are known for $A$ and $B$.
They are, most of the time, not already signed distance functions.
Note that the grid was chosen such that the origin (and thus the optimal $x$
for small displacements) is at the \emph{centre} of a grid cell and can
not be resolved exactly by the grid.
This yields the ``plateau'' in the error for small displacements.
However, as soon as external Hausdorff distances are admitted, the observed
error falls rapidly in accordance to \thref{ext dH estimate}.
The ``steps'' in the blue and red lines are caused by the discrete
nature of the grid.
The black curve shows the expected upper bound, which is given
by $\Delta_2$ for small displacements and by \eqref{ext dH estimate}
for larger ones.
(For our example situation, the maximum allowed radius $r$
in \defref{ext dH} can be computed exactly.)
One can clearly see that the theoretical and numerical results match
very well in their qualitative behaviour.

\includefig{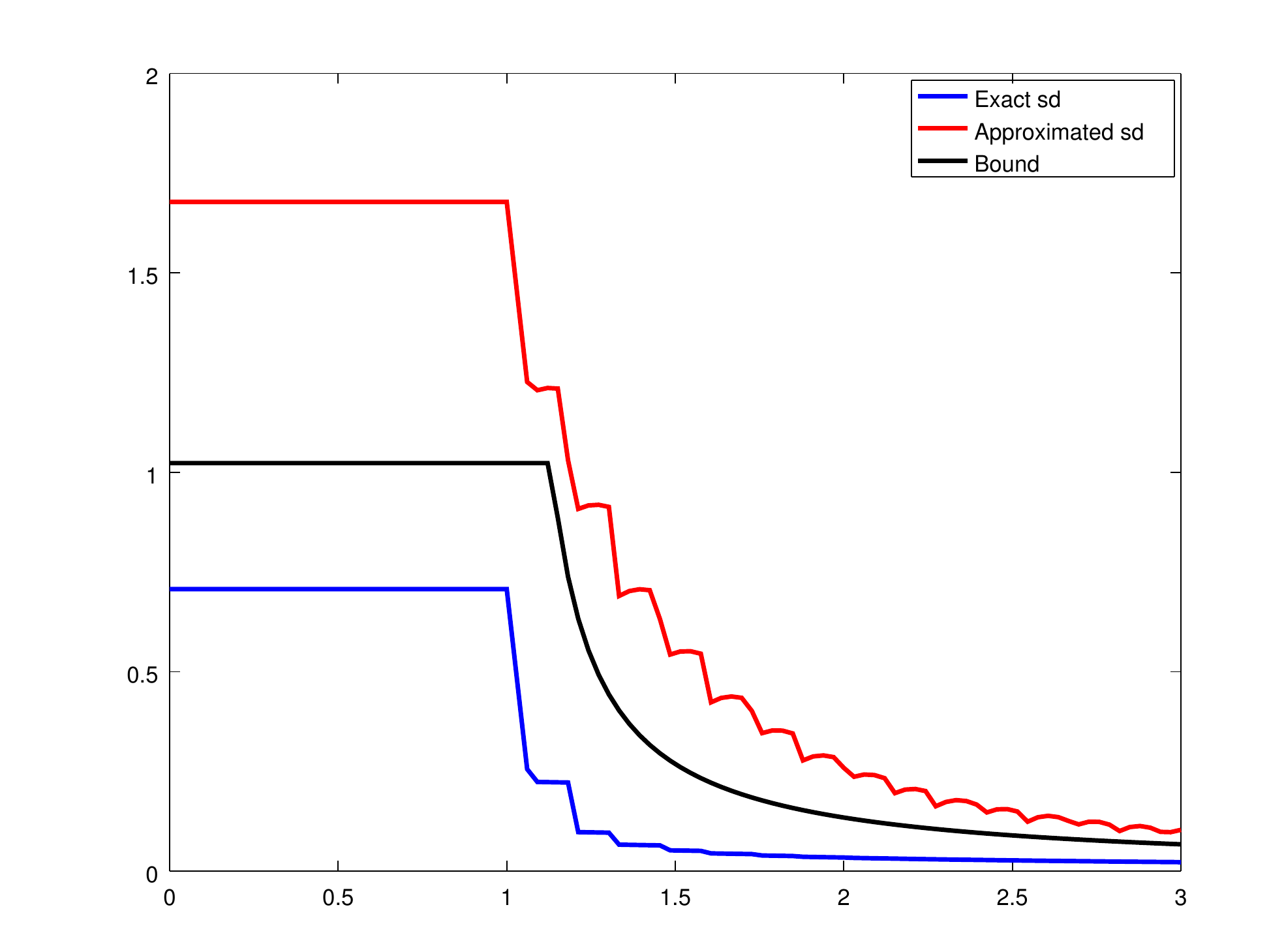}{width=0.8\textwidth}{plot movement}
  {Error $\delta$ (relative to $h$) for a situation similar to
   \figref{circle in ring}.
   The blue curve is based on exact signed distance functions.
   For the red curve, $\sd{A}$ and $\sd{B}$ were approximated as well.
   Black shows a combination of the error estimates
   from \corref{bound for suitable grid} (small displacement) and
   \eqref{ext dH estimate} (where applicable for larger displacements).
   The $x$-axis shows how far the inner circle is moved away from the origin.}

\figref{plot h} shows how the error depends on the grid spacing $h$.
The blue and red data is as before.
In the upper \subfigref{plot h}{centre}, the inner circle is
at the origin.  This is the situation of \subfigref{circle in ring}{centre},
and corresponds to the very left of the curves in \figref{plot movement}.
Here, the upper bound of \corref{bound for suitable grid}
applies and is shown with the black curve.
The convergence rate corresponds to $\bigO{h}$,
which can be seen clearly in the plot.
On the other hand, the lower \subfigref{plot h}{moved} shows
how $\delta$ behaves if we have an external Hausdorff distance.
It corresponds to the very right in \figref{plot movement},
with the inner circle displaced from the origin similarly to
\subfigref{circle in ring}{moved}.
Here, \thref{ext dH estimate} implies $\bigO{h^2}$ convergence.
Also this is, indeed, confirmed nicely by the numerical calculations.
(While both plots look similar, note the difference in the scaling
of the $y$-axes!)

\includesubfig{plot h}
  {Dependence of the error $\delta$ on the grid spacing $h$.
   The three data series are as in \figref{plot movement}.}
  {
    \subfig{0.8}{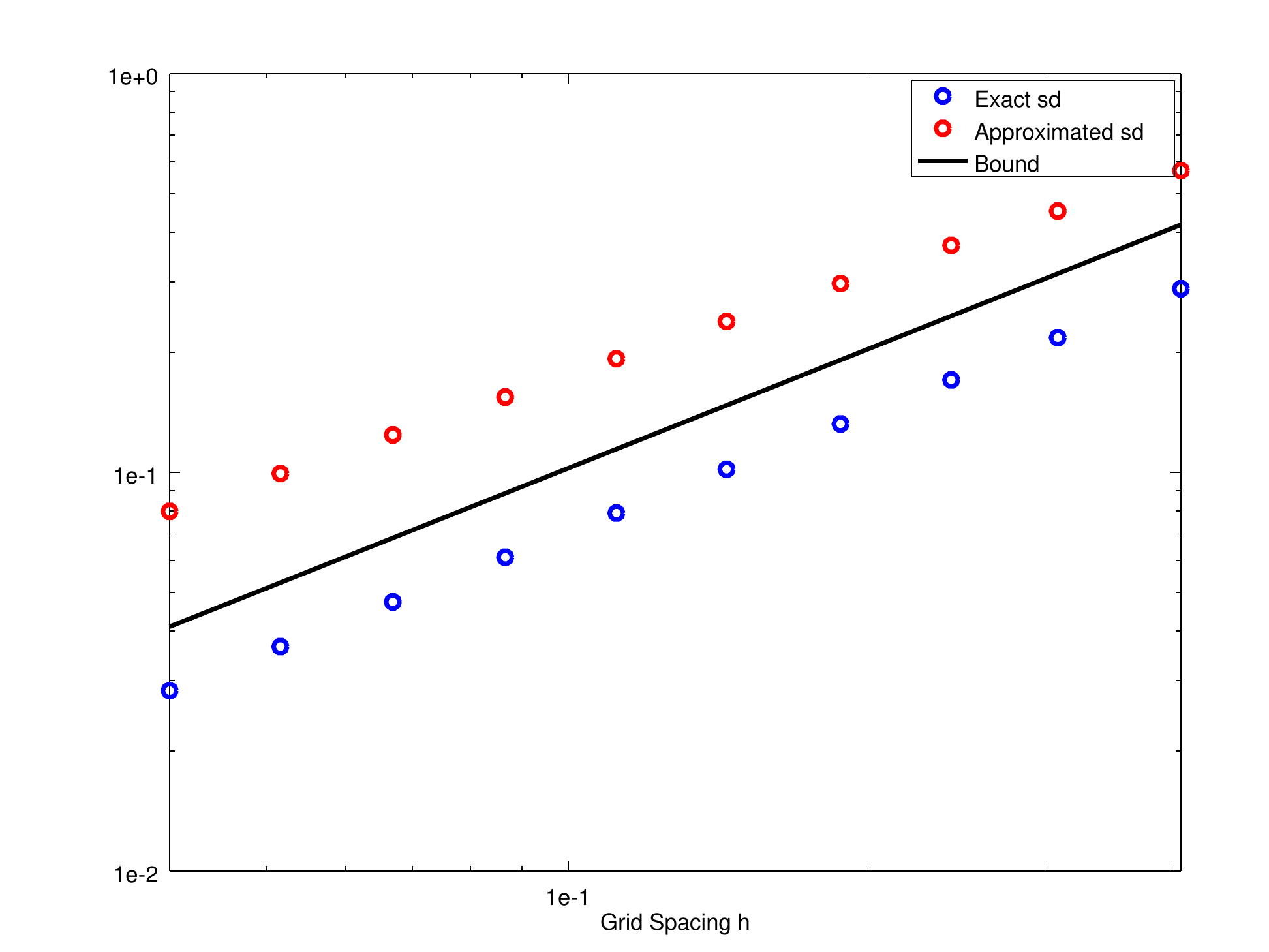}{centre}
      {Inner circle at the centre, as in \subfigref{circle in ring}{centre}.} \\
    \subfig{0.8}{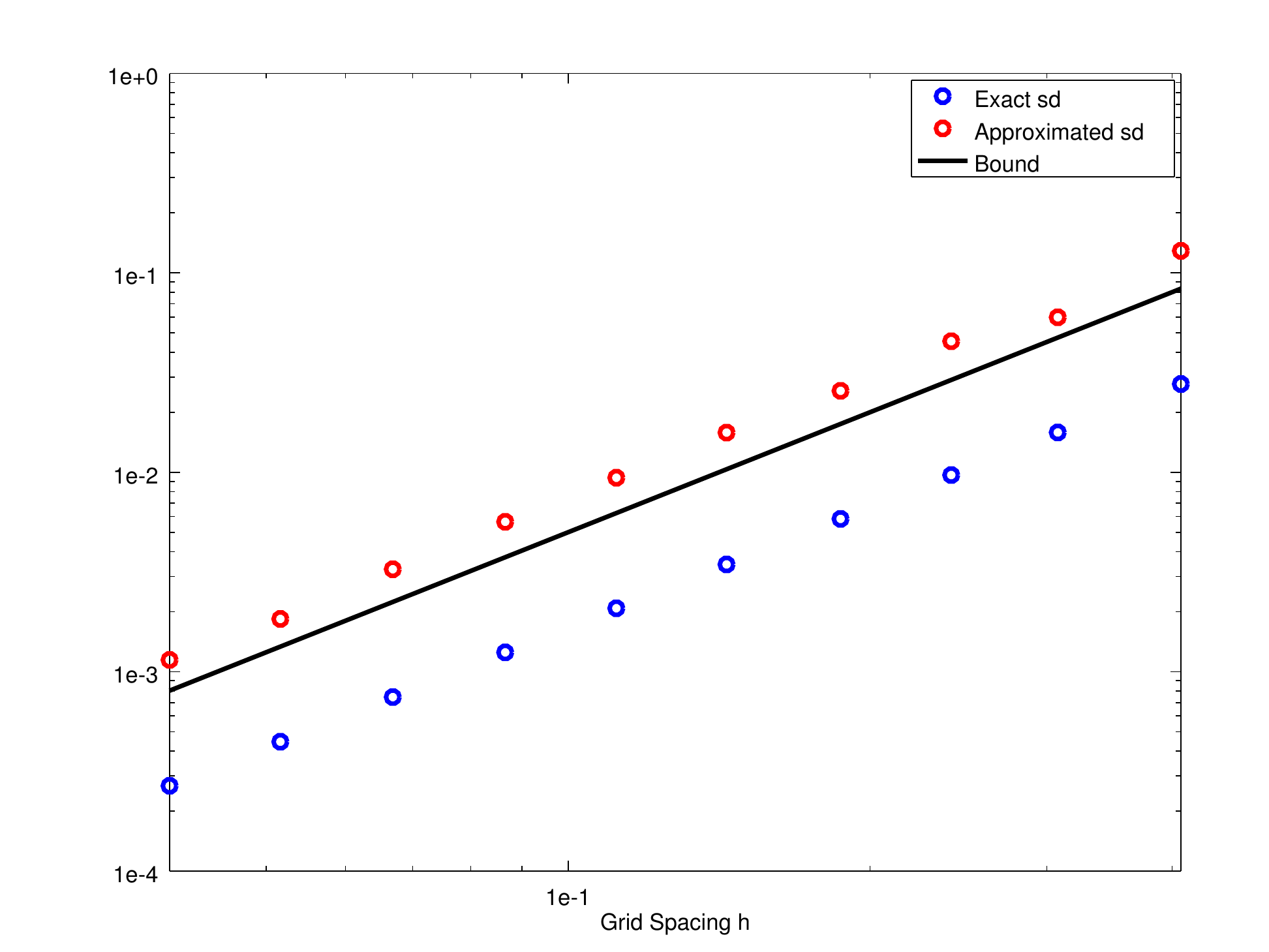}{moved}
      {Displaced inner circle, like \subfigref{circle in ring}{moved}.} \\
  }

\mysection{randomisation}{Improvements by Randomising the Grid}

Let us now take a closer look at the concept of external Hausdorff distances
and, in particular, the error estimate in the proof of \thref{ext dH estimate}.
An important ingredient for the resulting estimate (and the actual error)
is \emph{how close grid points come to lie to the external line}.
We can emphasise this even more by reformulating the error estimate
in the following way:


\begin{lemma}{dH estimate min grid point distance}
Let $A$ and $B$ admit an external Hausdorff distance with $r > 0$.
Choose $x \in A$, $y \in B$ and $d = (x - y) / \abs{x - y}$
as in \defref{ext dH}.
We denote by $\beta > 0$ the minimum distance any grid point
has to the part $L$ of the external line between $x + 3 / 8 \cdot r d$
and $x + 5 / 8 \cdot r d$, i.~e.,
\begin{equation}
\label{eq:dH estimate min grid point distance beta}
L = \set{x + \tau \cdot r d}{\tau \in \left[\frac38, \frac58\right]}, \;\;
\beta = \min_{z \in N} \d{L}(z)
  = \min_{z \in N} \min_{y \in L} \abs{y - z}.
\end{equation}
If the grid is large enough, then the error estimate
\begin{equation*}
\delta \le \frac{3}r \cdot \beta^2
\end{equation*}
holds for all grid spacings $h$ that are small enough.

\begin{proof}
We base the proof on \eqref{ext dH estimate delta}.
Using the notation of \figref{ext dH estimate}, let us
consider points $c'$ on the \emph{middle third} $M$ of the
external line.  For them, $\min(a, \, r-a) \ge r / 3$.
Consequently, \eqref{ext dH estimate delta} implies
\begin{equation*}
\delta \le \sqrt{b^2 + a^2} - a + \sqrt{b^2 + (r - a)^2} - (r - a)
  \le 2 \cdot \left(\sqrt{b^2 + \frac{r^2}9} - \frac{r}3\right)
  \le \frac3r \cdot b^2.
\end{equation*}
The last estimate can be seen with a series expansion.
This holds for $b$ being any distance (in normal direction)
of a grid point $z$ to $M$.
Furthermore, note that the interval $[3/8, \, 5/8]$ is a strict subset
of $[1/3, \, 2/3]$, which implies that also $L$ is a strict subset of $M$.
As long as $h$ is small enough, this means that we can find a suitable
grid point such that $b \le \beta$ holds.
This finishes the proof.
\end{proof}
\end{lemma}


Of course, the distance $\beta$ between grid points and the line segment $L$
can always be estimated trivially by $\sqrt{n} / 2 \cdot h$.  We did this
in the proof of \thref{ext dH estimate}, and this leads
precisely to the upper bound given in \eqref{ext dH estimate}.
The only difference to \lemmaref{dH estimate min grid point distance}
is the smaller constant in \eqref{ext dH estimate}.  This is due
to the very generous estimation
we used in \lemmaref{dH estimate min grid point distance} for simplicity.
We can now make an important observation:
This is the absolutely worst-case estimate, which matches
a situation as shown in \subfigref{ext dH}{grid}.
If the external line is not running ``parallel'' to the grid,
we can expect that some of the grid points lie much closer to it
(see \figref{line in grid}).
Such a situation occurs, in particular, almost surely if the grid
is placed ``randomly'' with respect to the geometry.  This usually happens,
for instance, if the original input data for a computation stems
from real-world measurements in one way or another.

\includefig{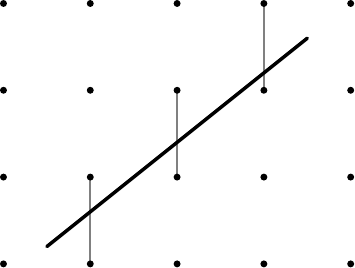}{width=0.35\textwidth}{line in grid}
  {Random placement of a line segment in a grid and the corresponding estimation
   of the distance $\beta$ to the closed grid point.}


In this section, we analyse the effect that such a random
grid placement has on the resulting error bounds.  We will see that
we can show even stronger convergence than in \thref{ext dH estimate}
in this case.  These estimates are based on
\lemmaref{dH estimate min grid point distance}.  The crucial additional
ingredient is a suitable upper bound for the minimum distance $\beta$.
The idea we employ for that is illustrated in \figref{line in grid}:
We look for edges of the grid that are intersected by the line segment $L$,
and use the distance \emph{along such an edge} to the next grid point
as an upper bound for $\beta$.  The finer the grid, the more such intersection
edges appear.  Each one of them gives us an additional ``chance'' to find
a particularly short distance, and thus improve $\beta$.
Consequently, it is interesting to know \emph{how many} such edges are there
for a particular line segment:

\begin{lemma}{intersected grid edges}
Let $L \subset \R^n$ be an arbitrary line segment with length $\lambda$,
placed in a grid as shown in \figref{line in grid}.
Then $L$ intersects at least
$\floor{\lambda / (\sqrt{n} \cdot h)}$
edges (or faces in 3D, cells in 4D and so on)
of the grid, where $h > 0$ is the grid spacing as usual.
In particular, the line $L$ of
\eqref{dH estimate min grid point distance beta}
intersects at least $r / (5 \sqrt n) \cdot 1 / h$ edges if $h$ is small enough.

\begin{proof}
Since $\sqrt{n} \cdot h$ is the longest diagonal of a grid cell,
every line segment of at least this length must intersect some edge.
Our line $L$ of length $\lambda$ is made up of
$\floor{\lambda / (\sqrt{n} \cdot h)}$ such pieces, which shows the claim.
The additional statement follows by noting that $\lambda = r / 4$
for $L$ in \eqref{dH estimate min grid point distance beta}
and that the estimate
\begin{equation*}
\floor{\frac{r}{4 \sqrt{n}} \cdot \frac1h}
  \ge \frac{r}{5 \sqrt{n}} \cdot \frac1h
\end{equation*}
is true if $h$ is sufficiently small.
\end{proof}
\end{lemma}


To further motivate this idea, let us continue
\subsecref{numerics} with more numerical examples.
We use the same basic situation (that is shown schematically
in \subfigref{circle in ring}{moved}), but add two randomised
changes:  First, the grid is randomly shifted in $x$- and $y$-direction
up to a grid cell.  This prevents any bias due to the placement
of the coordinate system relative to the outer circle.
Second, the inner circle is displaced in a \emph{random} direction,
which has the same effect as a random rotation of the grid.
The resulting errors $\delta$ between the exact and approximated
Hausdorff distances are plotted in \figref{random dir desc}.
These figures show the results of \num{1000} runs with different
randomisations.
The bound of \thref{ext dH estimate} is shown again with
the black line, and it clearly holds true for all runs.
Another observation, however, is that the error really decreases
much faster \emph{on average}:  The red dots indicate the
geometric mean over all runs; while not rigorously justified,
this gives a rough first indication of the average convergence behaviour.
A more representative analysis can be based on \figref{random dir expon}:
There, we show histograms of the approximated
convergence orders for each of the individual runs.
One can clearly see that most of the runs have, indeed, a stronger
order than $\bigO{h^2}$.  For 2D, the median order is almost $\bigO{h^4}$.
For 3D, it is ``only'' about $\bigO{h^3}$, but that is still a drastic
improvement over \thref{ext dH estimate}.
Also note that a tiny fraction of runs shows an order \emph{worse}
than $\bigO{h^2}$.  This is, however, not in contradiction
to our theory and rather a numerical artefact:  For them, the error is
(due to good luck) already very small for coarse grids,
so that the corresponding decrease seems slower for the grid sizes
we considered.  The error bound of \eqref{ext dH estimate}
is, nevertheless, still true.

\includesubfig{random dir desc}
  {The errors $\delta$ between the exact and approximated Hausdorff
   distances for the situation of \subfigref{circle in ring}{moved}
   with randomisation applied.  The black line shows the bound
   of \thref{ext dH estimate}, the dots correspond to the results
   of individual runs.  Their geometric mean is shown with red circles.}
  {
    \subfig{0.8}{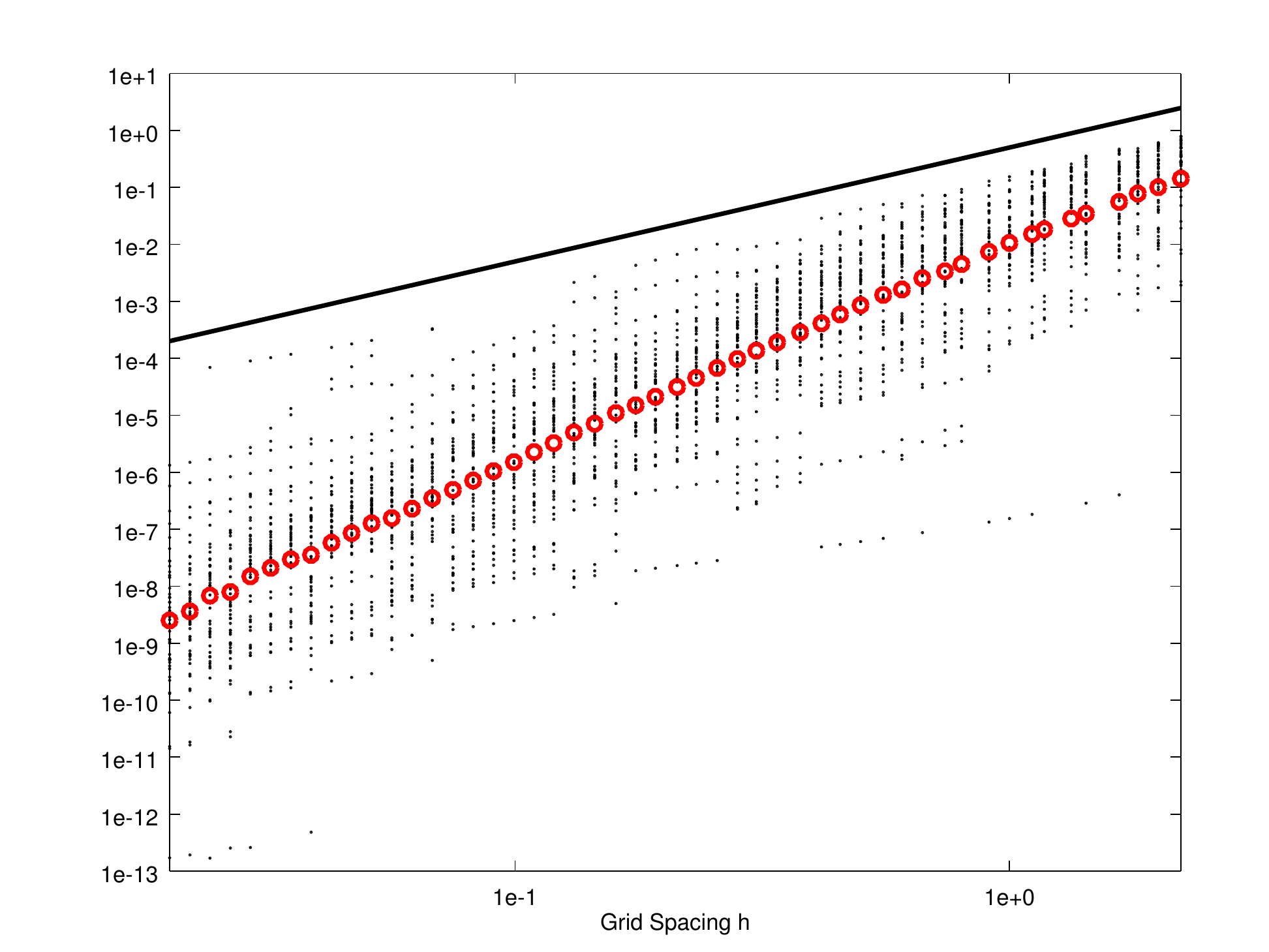}{2D}{2D} \\
    \subfig{0.8}{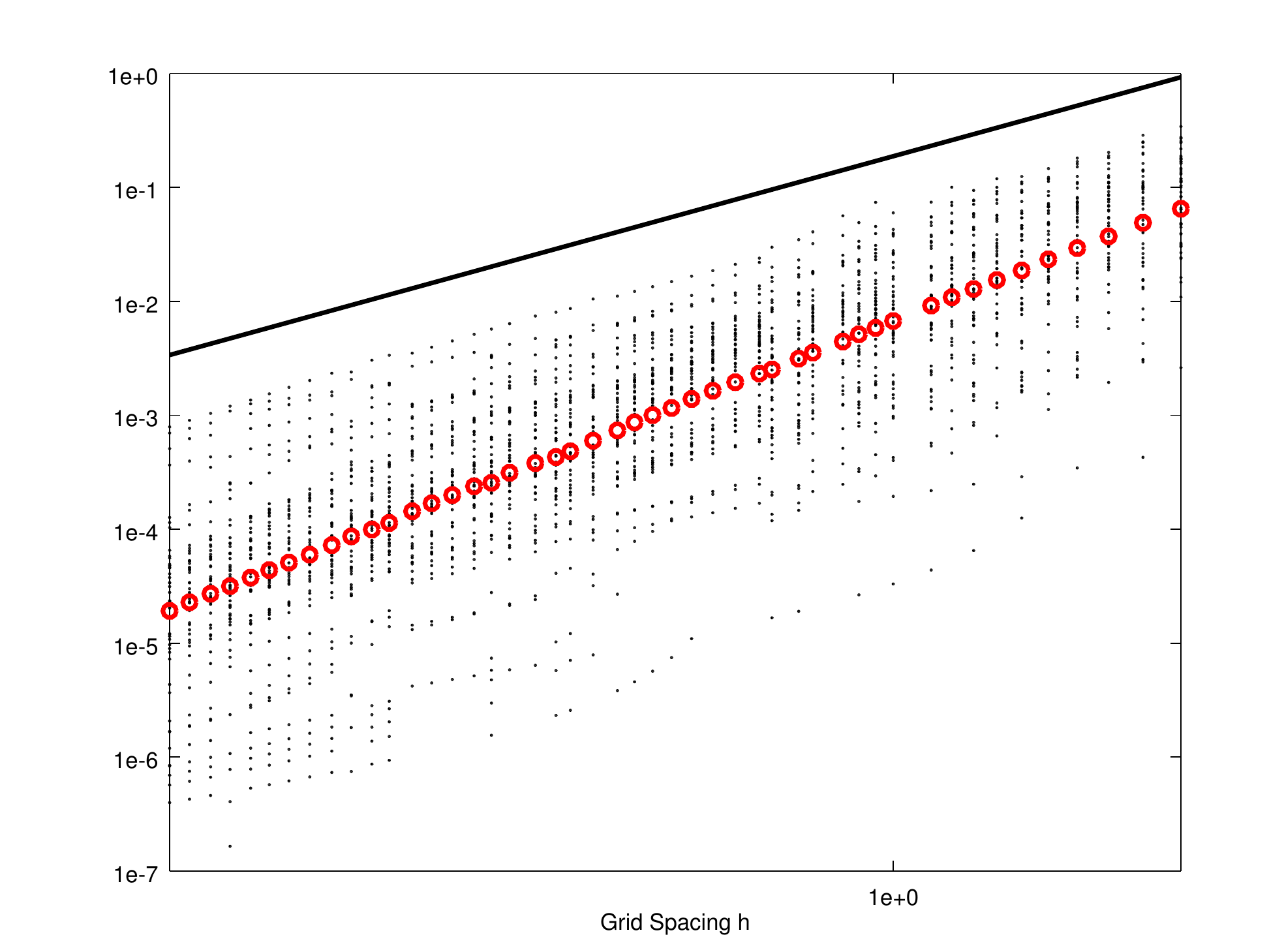}{3D}{3D} \\
  }

\includesubfig{random dir expon}
  {Histogram of the approximated convergence orders for the
   runs in \figref{random dir desc}.}
  {
    \subfig{0.8}{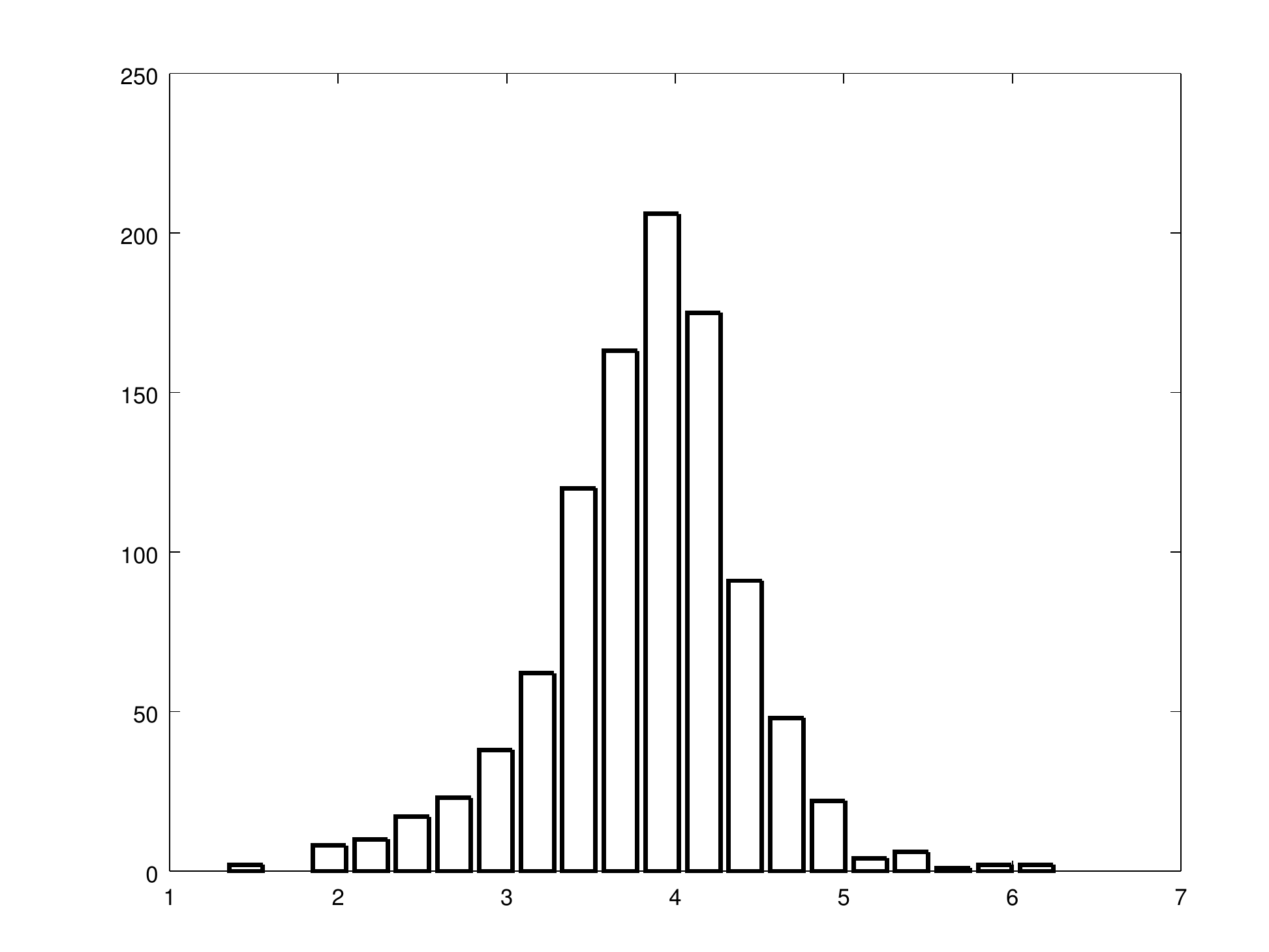}{2D}{2D} \\
    \subfig{0.8}{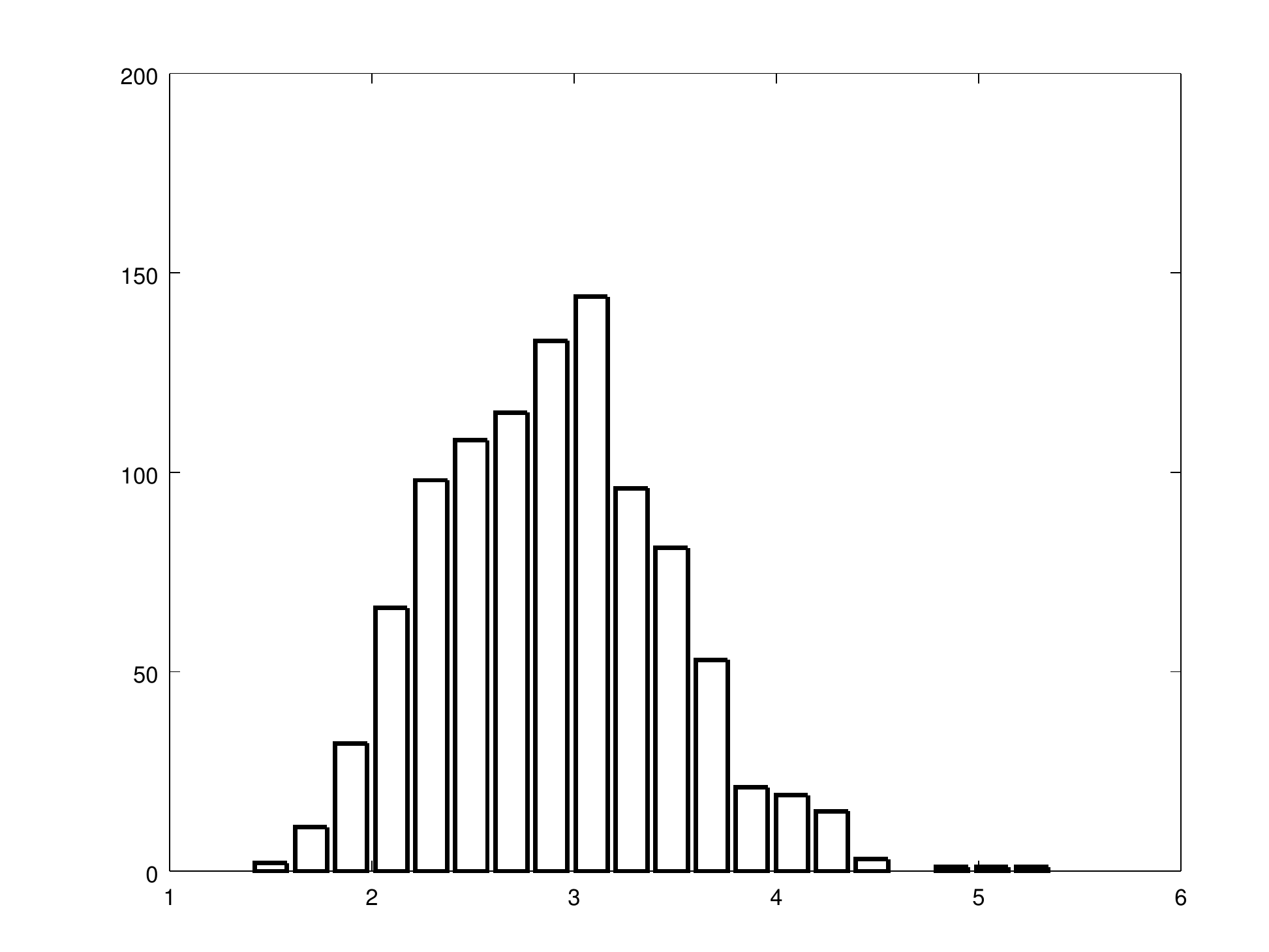}{3D}{3D} \\
  }

\mysubsection{2d almost-sure}{Super-Quadratic Convergence in 2D}

Let us now consider the randomised situation more thoroughly.
In this subsection, we concentrate on the 2D case and will
be able to show improved error bounds that hold
almost surely if a randomisation is applied.
These results will be complemented by a heuristic argument
given in \subsecref{unf distribution} that yields even better
convergence rates for the ``average case'' and works
in arbitrary dimension.


The main thing to do now is to find a bound
on $\beta$ of \eqref{dH estimate min grid point distance beta}.
For this, consider \figref{line in grid} again and
assume $h = 1$ for the moment.
The distance between the line segment $L$ and the grid point
below it along each intersected edge behaves then like
$\floor{x_0 + i k}$, where $x_0, k \in \R$ depend on the particular
line in question and $i \in \N$ numbers the intersected grid edges.
The analysis of such a sequence of numbers lies at the heart
of the error bound shown later in this subsection:

\begin{lemma}{d + i k minimum}
Let $x_0 \in \R$ and $k \in \R \setminus \Q$ be given.
For $i \in \N$, we define $x_i = \floor{x_0 + i k} \in [0, 1)$.
Then, for arbitrary $N \in \N$, there are
$i_0, j_0 \in \directset{0, \ldots, N}$
such that $0 < \epsilon = \abs{x_{i_0} - x_{j_0}} \le 1 / N$.
Furthermore, there is $m \in \directset{0, \ldots, K}$
with $\abs{x_m} \le \epsilon$, where the maximum number $K$
of iterates necessary for achieving this condition satisfies
\begin{equation*}
K \le N \ceil{\frac1\epsilon} \le \frac2{\epsilon^2}.
\end{equation*}

\begin{proof}
Assume for a moment that $\abs{x_i - x_j} > 1 / N$ for all $\intrange{i,j}0N$
with $i \ne j$.
This also implies that all intervals $I_i = [x_i - 1 / N, \, x_i + 1 / N]$
are disjoint.  Since $\abs{I_i \cap (0, 1)} \ge 1 / N$
and there are $N + 1$ of the intervals, this is not possible.
Thus, there exist $i_0, j_0 \in \directset{0, \ldots, N}$
with $i_0 \ne j_0$ and $\epsilon = \abs{x_{i_0} - x_{j_0}} \le 1 / N$.
It remains to show that $\epsilon > 0$.
Assume to the contrary that $\abs{x_{i_0} - x_{j_0}} = 0$.
This means $\abs{i_0 - j_0} k \in \N$.
But since $i_0 \ne j_0$, this contradicts the assumption that $k \not\in \Q$.
Thus, the first statement is shown.

Without loss of generality, let us now assume that
$i_0 < j_0$ and that $x_{i_0} < x_{j_0}$.
Setting $p = j_0 - i_0 \in \N$, note that
$p \in \directset{1, \ldots, N}$ and
$\floor{p k} = \abs{x_{i_0} - x_{j_0}} = \epsilon$.
In other words, advancing $p$ iterates in the sequence increases
the value of $x_i$ by $\epsilon$.
With a similar argument to before, this implies
that we cover the whole interval $[0, 1]$ in at most
$p \ceil{1 / \epsilon}$ iterations.
This implies that $\abs{x_m} \le \epsilon$
for some $m \in \directset{0, \ldots, K}$ with
\begin{equation*}
K \le p \ceil{\frac1\epsilon}
  \le N \ceil{\frac1\epsilon}
  \le \frac1\epsilon \left(\frac1\epsilon + 1\right)
  \le \frac2{\epsilon^2}.
\end{equation*}
\end{proof}
\end{lemma}

There are a few things to remark about the proof
of \lemmaref{d + i k minimum}:
Note that we do not directly get an upper bound for the number
of iterates required to get close to zero within a certain
threshold.  Instead, the sequence itself yields a distance $\epsilon$,
and the number of iterates required is then defined in terms of $\epsilon$.
If $\epsilon$ happens to be much smaller than $1 / N$, we also need a lot
more iterates.  On the other hand, however, this larger number
of iterates also brings us much closer to zero than $1 / N$.
There is a balance between closeness to zero and the number
of iterates we need, but we have no direct control over either quantity.
This also shows why it is crucial that $k$ is an \emph{irrational} number:
If it is not, it may happen that $\epsilon = 0$.  In this case,
we are not guaranteed to ever come close to zero in a finite number
of iterations since the sequence becomes periodic.
In terms of our error analysis, this corresponds to the situation
that the external line is parallel to the grid as depicted,
for instance, in \subfigref{ext dH}{grid}.
Fortunately for us, however, the rational numbers have measure zero.
This means that a randomly rotated grid will almost surely yield
an irrational $k$ so that our analysis applies.

Another interesting observation is the following:
The important first part of the proof establishes an estimate on
the number of steps necessary before a pair of iterates occurs that
are close to each other.  On the first thought, this may sound as if
the \emph{birthday paradox} is applicable in this situation.  This would
improve the number of iterates necessary to $\bigO{\sqrt{N}}$.
Unfortunately for us, however, this is not true:
In our case, the actual question is
\emph{when $\floor{p k}$ comes close to zero}.
It is only formulated in terms of a \emph{pair} of iterates
because this seems like the more natural formulation for the proof.

Finally, let us briefly discuss what happens if we try to use
the same approach for the Hausdorff distance in $\R^3$:
In this case, we consider
a sequence of points on the unit square $[0, 1)^2$.
With the same argument as in the proof of \lemmaref{d + i k minimum},
we can still show that a pair of iterates is within $1 / N$ of each other
if we perform $\bigO{N^2}$ steps.  (Note that the estimate gets drastically
worse here!)  It is not so clear, however, whether the second part
of the proof can be adapted:  On the one-dimensional interval $[0, 1)$,
repeatedly adding $\epsilon$ to an arbitrary point will eventually wrap around
and yield an iterate close to zero.
On the square in 2D, however, there is much more freedom for the sequence
to iterate without ever getting close to the origin.


With the technical result of \lemmaref{d + i k minimum}
in place, we can now use it to derive an estimate for $\beta$.
Combining this with \lemmaref{dH estimate min grid point distance}, we obtain
an improved error estimate $\delta$ for the approximate Hausdorff distance:

\begin{corollary}{superquadratic estimate}
Let $n = 2$ and $A, B \subset \R^2$ admit an external Hausdorff distance
with $r > 0$.  Assume that the grid is rotated randomly
with respect to the external line.
Then, with probability one, there exists a sequence $\seqk{m}{h_m} \to 0$ of
grid spacings tending to zero such that
\begin{equation}
\label{eq:superquadratic estimate}
\delta \le \frac{45}{r^2} \cdot h_m h^2
\end{equation}
for all $m \in \N$ and grids with $h < h_m$.
In particular, the approximation error $\delta$
vanishes \emph{super-quadratically} in the grid spacing.

\begin{proof}
For $m \in \N$, choose $\epsilon_m \in (0, \, 1 / m]$
according to \lemmaref{d + i k minimum} and define
\begin{equation*}
h_m = \frac{\epsilon^2 r}{10 \sqrt{n}}
  \;\; \Leftrightarrow \;\;
  \epsilon_m = \sqrt{\frac{10 \sqrt{n}}r} \cdot \sqrt{h_m}.
\end{equation*}
Then clearly $h_m \to 0$ as $m \to \infty$.
Furthermore, this definition ensures that
the line segment $L$ of \eqref{dH estimate min grid point distance beta}
intersects at least $2 / \epsilon_m^2$ edges
of any grid with $h \le h_m$ according to
\lemmaref{intersected grid edges}.

As discussed above, the intersections of the line segment $L$ with grid edges
(see \figref{line in grid}) can be described by a sequence
like $\floor{x_0 + i k}$ after scaling the grid cells to size one.
Note that a random rotation of the grid ensures $k \in \R \setminus \Q$
with probability one.
Thus, \lemmaref{d + i k minimum} is applicable and yields that
$\beta \le \epsilon_m h$, where the additional factor $h$
can be added since our grid has cells of size $h \times h$
instead of unit size.
Thus, \lemmaref{dH estimate min grid point distance} yields
\begin{equation*}
\delta \le \frac3r \cdot \epsilon_m^2 h^2
  = \frac{30 \sqrt{2}}{r^2} \cdot h_m h^2
  \le \frac{45}{r^2} \cdot h_m h^2.
\end{equation*}
\end{proof}
\end{corollary}

Let us emphasise again that \corref{superquadratic estimate} only shows
\emph{super-quadratic} convergence and
$\bigO{h^3}$ along a \emph{particular sequence} of grids, not
cubic convergence in general.
Since \lemmaref{d + i k minimum} does not give us any control
over $\epsilon$, we also do not have any knowledge about
the resulting sequence $\seqk{m}{h_m}$ along which third-order
convergence occurs.  In practice, however, the bound
in \eqref{superquadratic estimate} seems to be quite conservative
and usually an overestimation of the error.
This matches also the behaviour seen in \subfigref{random dir expon}{2D},
which suggests that most of the runs show an empirical convergence order
stronger than $\bigO{h^3}$.

\mysubsection{unf distribution}{A Heuristic, Statistically-Motivated Estimate}

While \corref{superquadratic estimate} gives a rigorously proven
estimate, it is not clear how to extend the result to situations
in more than two dimensions.  Furthermore, the result is rather
a conservative worst-case estimate than a practical bound
on the expected average error.  In this subsection, we give
a different argument leading to stronger
estimates that match the empirical results of \figref{random dir expon}
more closely.
This argument works in arbitrary space dimension, although the statement
gets weaker the higher the dimension is.
Also note that the main idea is only heuristically motivated
and cannot be recovered in a formal proof.
Consequently, all of this subsection will be written in an informal
way without rigorously proving any statements.

We still base our estimate on \lemmaref{dH estimate min grid point distance},
but replace \lemmaref{d + i k minimum} by the following heuristic idea:
If the line segment $L$ intersects a large enough number of grid edges
and is not parallel to the grid in any way, it seems plausible to assume
that the intersection distances $\floor{x_0 + i k}$
follow roughly a \emph{uniform distribution} on $[0, 1)$.
In a rigorous setting, this is of course not true---we do not consider
independent random numbers at all, but correlated iterates of a
deterministic sequence.  Empirically, however, this assumption
seems to be justified relatively well:
The sequence can be compared to the well-known
\emph{linear congruential generators}
for pseudo-random numbers (see, for instance,
Section~3.1 of \cite{sheldonSimulation}).
They are among the simplest and statistically
weakest PRNGs, but their weaknesses seem to be irrelevant for our situation.
See also \figref{unf dist}:  This plot shows a histogram of \num{10000}
iterates of such a sequence.  One can clearly see that while the distribution
seems not to be perfectly uniform, it is quite close.  The more iterates
we use, the more uniform the distribution gets.
Since we are interested in the limit of fine grids, i.~e., $h \to 0$,
this corresponds to a large number of intersected edges.
Thus, it seems reasonable to assume that the intersections between
the line segment $L$ and grid edges follow such a uniform distribution.

\includefig{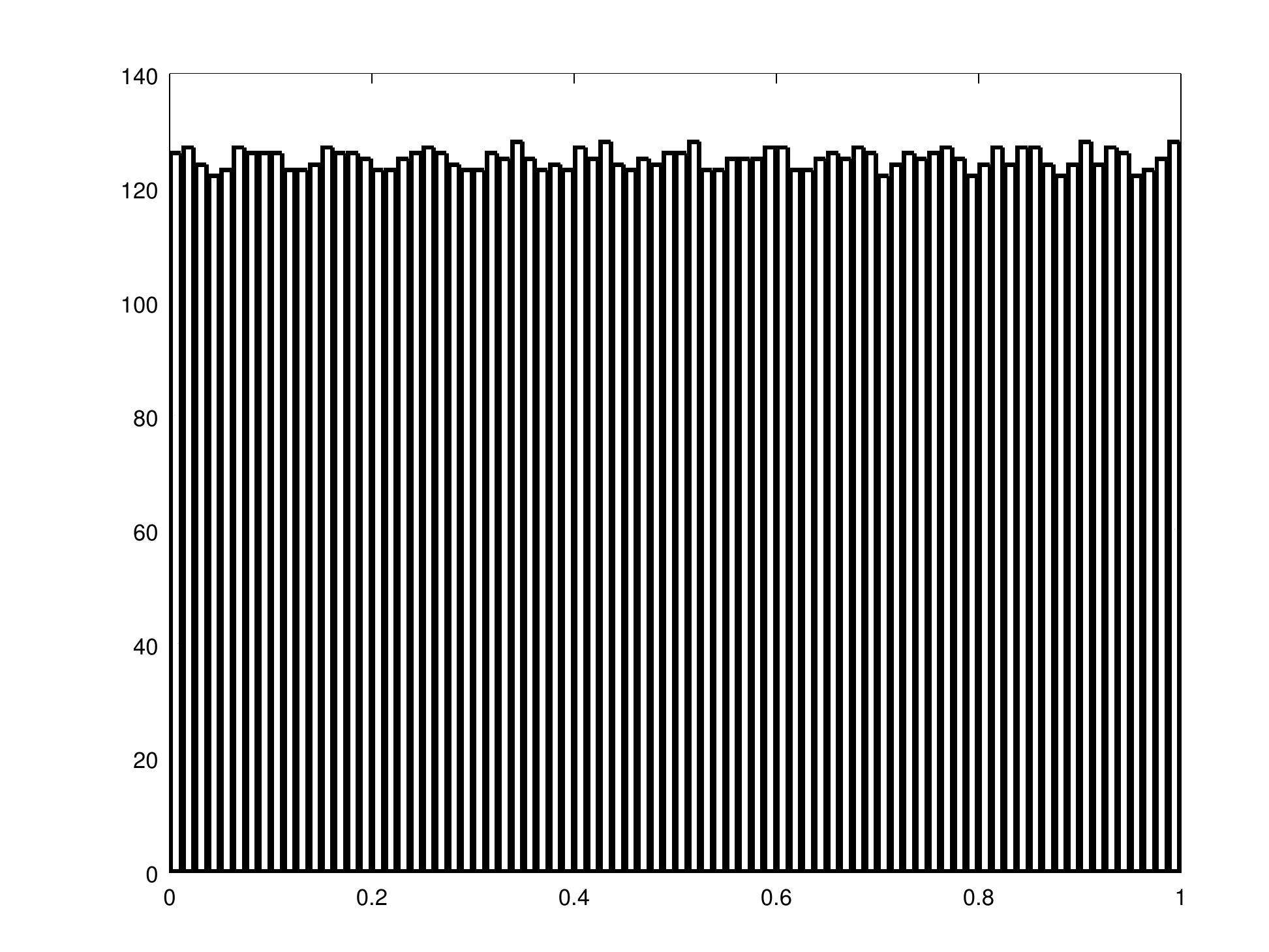}{width=0.8\textwidth}{unf dist}
  {Histogram of \num{10000} iterates in a sequence
   $\floor{x_0 + i k}$ with randomly chosen $x_0, k \in [0, 1]$.}

This idea can be applied also in three and more dimensions:
Instead of intersections along grid edges, we consider
intersections of $L$ with faces of the grid cells in 3D.
For them, we may (with the same heuristic justification) assume that
they follow a uniform distribution on the square $[0, 1)^2$.
In even higher dimensions, this approach can be generalised further.
Thus, let us assume that $X_1, \ldots, X_N$ is a sequence
of independent random variables that follow a uniform distribution
on $[0, 1)^{n-1}$.
To replace \lemmaref{d + i k minimum} and find an estimate on $\beta$
for \lemmaref{dH estimate min grid point distance}, we have to analyse
the statistical properties of their minimum distance to the origin,
i.~e., of the derived random variable
\begin{equation*}
X = \min(\abs{X_1}, \, \ldots, \, \abs{X_N}).
\end{equation*}
For simplicity, we will not analyse $X$ directly.
Instead, we base our computations on random variables
$Y_i$ that are uniformly distributed on the set
\begin{equation}
\label{eq:unf dist sector Q_n}
Q_n = \set{x \in [0, \infty)^{n-1}}{\abs{x} < \sqrt{n - 1}} \subset \R^{n-1}.
\end{equation}
This is a superset of $[0, 1)^{n-1}$ and a spherical sector.
For 3D, the set $Q_3$ is illustrated in \figref{unf dist sector}.
Since the additional regions of $Q_n$ with respect
to $[0, 1)^{n-1}$ are far away from the origin, the corresponding
minimum $Y$ is ``larger than'' $X$ with respect to important properties
such as the expectation value or quantiles.
Furthermore, with an appropriate scaling, one can also
produce a \emph{subset} of $[0, 1)^{n-1}$ from $Q_n$.
This implies that $X$ and $Y$ differ, roughly speaking,
at most by a dimensional constant.
Thus, we can use $Y$ for our analysis of the resulting convergence order.

\includefig{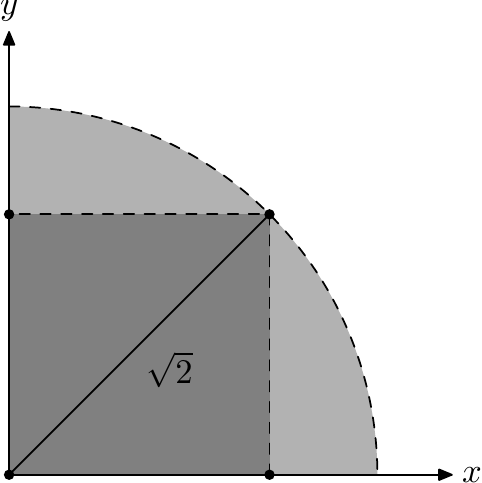}{width=0.3\textwidth}{unf dist sector}
  {Relation between the quarter circle $Q_3$ of \eqref{unf dist sector Q_n}
   and the square $[0, 1)^2$.  The latter is filled
   dark grey, with the additional regions of $Q_3$ light grey.}

Next, we compute the \emph{density function} for a random variable
$\abs{Y_i}$.  For this, we are interested in the question
what distance a uniformly placed point on $Q_n$ has from the origin.
For $n \in \N$, let us denote the volume of the $n$-dimensional unit ball
$\B10 \subset \R^n$ by $\omega_n$.  For a general derivation of these
constants, see Theorem~26.13 in \cite{yehMeasure}.
Furthermore, note that the \emph{surface measure} of a spherical shell
with radius one is given by $n \omega_n$ according
to Observation~26.24 in \cite{yehMeasure}.
It is then straight-forward to express the density function $f(r)$
for a particular distance $r = \abs{Y_i}$ with these constants as
\begin{equation*}
f(r) = \frac{r^{n-2} \cdot (n - 1) \omega_{n-1}}
            {\left(\sqrt{n - 1}\right)^{n-1} \cdot \omega_{n - 1}}
  = (n - 1)^{\frac{3-n}2} \cdot r^{n-2}.
\end{equation*}
Note that it is trivial to check that this density function
is normalised, i.~e.,
\begin{equation*}
\integral0{\sqrt{n-1}}{f(r)}{dr}
  = (n-1)^{\frac{3-n}2} \integral0{\sqrt{n-1}}{r^{n-2}}{dr} = 1.
\end{equation*}
Building on $\abs{Y_i}$, it remains to find the distribution
of the \emph{minimum} $Y$ of $N$ such random variables.
In the case $Y = \abs{Y_1}$, we know that $Y$ takes the norm of $Y_1$
and that all other points $Y_2, \ldots, Y_N$ must have a norm at
least as large.  The probability density for this to happen
is consequently
\begin{equation*}
\tilde{g}(r) = f(r) \left(\integral{r}{\sqrt{n-1}}{f(\rho)}{d\rho}\right)^{N-1}
  = (n-1)^{N \frac{3-n}2} \cdot r^{n-2}
      \cdot \left((n-1)^{\frac{n-3}2} - \frac{r^{n-1}}{n-1}\right)^{N-1}.
\end{equation*}
Since any of the $N$ variables could be the minimum,
the density function of $Y$ is given by
\begin{equation*}
g(r) = N \cdot \tilde{g}(r)
  = N \cdot f(r) \left(\integral{r}{\sqrt{n-1}}{f(\rho)}{d\rho}\right)^{N-1}.
\end{equation*}
Again, one can show that $g$ is normalised by integrating
over $r \in [0, \sqrt{n-1})$.

Let us now compute the \emph{expectation value} $\exptval{Y}$.
This quantity tells us how close the minimum distance of a grid
point is \emph{on average} to the line segment $L$.
Again, the computation is straight-forward but technical.
One finds
\begin{equation*}
\exptval{Y}
  = \integral0{\sqrt{n-1}}{r \cdot g(r)}{dr}
  = N \sqrt{n-1} \cdot \Beta{\frac{n}{n-1}}N
  = N \sqrt{n-1} \cdot \gammafcn{\frac{n}{n-1}}
      \cdot \frac{\gammafcn{N}}{\gammafcn{N + \frac{n}{n-1}}},
\end{equation*}
where $\Beta\cdot\cdot$ and $\gammafcn\cdot$ are the beta
and gamma functions, respectively.
See Chapter~5 of \cite{dlmf} for more details about these
special functions.
We are particularly interested in the limit of fine grids,
corresponding to $N \to \infty$.
The asymptotic behaviour in this limit can be derived
from 5.11.12 in \cite{dlmf}, yielding
\begin{equation}
\label{eq:exptval Y}
\exptval{Y}
  \sim N \cdot \frac{\gammafcn{N}}{\gammafcn{N + \frac{n}{n-1}}}
  \sim N \cdot N^{-\frac{n}{n-1}}
  \sim N^{-\frac1{n-1}}.
\end{equation}
For $n = 2$, one can also compute the percentiles of $Y$
exactly.  They show precisely the same asymptotic behaviour.

Finally, let us discuss what \eqref{exptval Y} implies
for the convergence order of $\delta$.
For this, note first that $N \sim 1 / h$ according to
\lemmaref{intersected grid edges}.
The minimum distance $\beta$ of
\eqref{dH estimate min grid point distance beta} is given
by $h X$ after scaling the grid cells to unit size
(matching the definition of $X$).
As discussed above, we approximate this value by
\begin{equation*}
\beta \approx h \cdot \exptval{Y}
  \sim h \cdot h^{\frac1{n-1}}.
\end{equation*}
Thus, \lemmaref{dH estimate min grid point distance}
implies that we can expect the approximation error to behave like
\begin{equation*}
\delta \sim \beta^2
  \sim h^2 \cdot h^{\frac2{n-1}}
    =  h^{\frac{2n}{n-1}}.
\end{equation*}
In other words, we get, indeed, a stronger average convergence
order than $\bigO{h^2}$ of \thref{ext dH estimate}.
The additional factor shrinks with increasing space dimension.
For $n = 2$, we get $\bigO{h^4}$.  With $n = 3$, the expected
order is $h^3$.
Note that this matches precisely the empirical results
found in \figref{random dir expon}.

\mysection{outlook}{Conclusion and Outlook}

In this paper, we have discussed the relation between the Hausdorff
distance $\dH(A, B)$ of two sets and their (signed) distance functions.
This allowed us to state a method for computing $\dH(A, B)$ based on
these distance functions $\d{A}$ and $\d{B}$.
They can be easily computed, for instance, if the sets themselves
are described in a level-set framework.
Furthermore, we were able to analyse the approximation error
made if the distance functions are only known on a finite grid
(which is mostly the case in applications).
To the best of our knowledge, such an error analysis has not been
carried out before.
We were able to derive a general and sharp estimate
in \corref{bound for suitable grid}.
For the more regular situation characterised in \defref{ext dH}, our result of
\thref{ext dH estimate} implies an even lower error bound:
In this case, the error is at most $\bigO{h^2}$, where $h$
is the grid spacing.
With a random rotation of the grid, even better convergence
rates can be achieved.  We have rigorously shown super-quadratic convergence
for geometries in $\R^2$ in this case, and given a heuristic derivation
of the average convergence order.
All of our results are confirmed by numerical experiments.
Nevertheless, there remain a lot of areas open for further investigations.
In particular, we believe that the following open questions
would be interesting targets of further research:

\begin{itemize}

\item
We have not given a general formula for the constants
$\Delta_n$ in arbitrary space dimensions.
While cases with $n > 3$ are, arguably, not so important,
it could still be interesting to consider a general dimension $n$.
It is probably not too hard to devise a method of assembling
the system of equations that characterises $\Delta_n$
similarly to \eqref{maximal error p}.

\item
Going even further, is it possible to find an efficient
algorithm to evaluate \eqref{general upper bound cell}
in a general situation?
If this was possible, one could exploit \thref{general upper bound}
directly instead of resorting to approximations like
\corref{bound for suitable grid}.
This is particularly interesting, since it would allow to
couple information about $\abs{\d{A}(x) - \d{B}(x)}$ on the corners $x$
of a grid cell with the corresponding Lipschitz constants in $t(x, \cdot)$.
Doing this on a per-cell basis could improve the error bounds even further
for an actual computation.

\item
Our definition of external Hausdorff distances in \defref{ext dH}
is straight-forward to verify for concrete situations and we argued
why we believe that a lot of practical shapes fall into this category.
It would, nevertheless, be interesting to further analyse the class of
geometries that admit such an external Hausdorff distance.

\item
While we believe that our discussion of the effect of randomisation
in \secref{randomisation} forms a sound basis both theoretically and
for practical applications, the stochastic error analysis started
there opens up a lot of possibilities for further refinement.

\end{itemize}

\section*{Acknowledgement}

The author would like to thank Bernhard Kiniger (Technical University
of Munich) for pointing out that the complementary Hausdorff distance
is worth a discussion in its own right.
The improvements of \secref{randomisation} were only possible
due to Michael Kerber from the Technical University of Graz, who
suggested that a randomisation of the grid could be an interesting
further topic for research.
This work is supported by the Austrian Science Fund (FWF) and
the International Research Training Group IGDK~1754.

\bibliographystyle{plain}
\bibliography{references}

\end{document}